\documentclass[11pt,a4paper]{article}

\title{Characterization of valid auxiliary functions for representations of extreme value distributions and their max-domains of attraction}

\RequirePackage[colorlinks,linkcolor=blue,citecolor=blue,urlcolor=blue]{hyperref}

\usepackage{exscale,amsthm,amssymb}
\usepackage[dvips]{graphicx}
\usepackage{color}
\usepackage{float}
\usepackage[english]{babel}
\usepackage{rawfonts,latexsym,lscape,epsfig,multirow,lscape,graphics,bm,bbm}
\usepackage{rotating,times,amsmath}

\usepackage{natbib}
\usepackage{etoolbox}

\def\oF{\overline{F}}
\def\oH{\overline{H}}
\def\oG{\overline{G}}

\def\Pg{\mathcal{P}_\gamma}
\def\P{\mathcal{P}}
\def\g{\gamma}

\def\AMF{\mathcal{AM}(F)}
\def\AVF{\mathcal{AV}(F)}

\def\I{\int_x^{x_E} \oF(u)\mathrm{d}u}
\def\II{\int_x^{x_E} \int_v^{x_E} \oF(u)\mathrm{d}u\, \mathrm{d}v}
\def\III{\int_x^{x_E} \int_w^{x_E} \int_v^{x_E} \oF(u)\mathrm{d}u\, \mathrm{d}v\, \mathrm{d}w}
\def\Is{\int_x^{x_E}\!\! \oF(u)\mathrm{d}u}
\def\IIs{\int_x^{x_E} \!\!\!\int_v^{x_E} \!\!\oF(u)\mathrm{d}u\, \mathrm{d}v}
\def\IIIs{\int_x^{x_E}\!\!\! \int_w^{x_E}\!\!\! \int_v^{x_E}\!\! \oF(u)\mathrm{d}u\, \mathrm{d}v\, \mathrm{d}w}
\newcommand{\dps}{\displaystyle}
\def\VR{$\textsf{VR}$ }
\def\VRs{$\textsf{VR}$ }
\def\VRp{$\textsf{VR}$}
\def\VRsp{$\textsf{VR}$}
\def\vMR{$\textsf{vMR}$ }
\def\vMRs{$\textsf{vMR}$ }
\def\vMRp{$\textsf{vMR}$}

\theoremstyle{plain}
\newtheorem{theorem}{Theorem}
\newtheorem{lemma}{Lemma}
\newtheorem{proposition}{Proposition}
\newtheorem{corollary}{Corollary}

\theoremstyle{remark}

\newtheorem{remark}{Remark}

\date{October 27, 2023}

\renewcommand{\baselinestretch}{1.2}

\setcitestyle{aysep={}}

\begin{document}

\author{ $\quad$\\[-2mm]
Miriam Isabel Seifert$^{\,1}$ \\[5mm]
}

\thispagestyle{empty}
\maketitle
\vspace{-12mm}

\begin{minipage}{0.045\textwidth}
$\quad$
\end{minipage}
\begin{minipage}{0.9\textwidth}
$\quad$\\
\begin{abstract}
{\small
\noindent In this paper we study two important representations for extreme value distributions and their max-domains of attraction (MDA), namely von Mises representation (\vMRp) and variation representation (\VRp), which are convenient ways to gain limit results. Both \VR and \vMR are defined via so-called auxiliary functions $\psi$. Up to now, however, the set of valid auxiliary functions for \vMR has neither been characterized completely nor separated from those for \VRp. We contribute to the current literature by introducing ``universal'' auxiliary functions which are valid for both \VR and \vMR representations for the entire MDA distribution families.
Then we identify exactly the sets of valid auxiliary functions for both \VR and \vMRp.
Moreover, we propose a method for finding appropriate auxiliary functions with analytically simple structure and provide them for several important distributions.\\

}
\end{abstract}
\end{minipage}

\footnotetext[1]{\small $\;$Ruhr University Bochum, Faculty of Management and Economics,\\
 $\phantom{abci}$ Universit\"{a}tsstra\ss e 150, 44801 Bochum, Germany.\\
 $\phantom{abci}$ Email: miriam.seifert@rub.de}

\newpage

\renewcommand{\baselinestretch}{1.41}\normalsize

\flushbottom

\section{Introduction}\label{sect:1}

The classical extreme value theory deals with the analysis of asymptotic distributions of sample extremes, i.e. of the maximum (or the minimum) of independently identically distributed random variables $X_1,\ldots, X_n$ with distribution function~$F$ for the sample size $n\to\infty$. Consider a distribution~$F$ for which,  under an appropriate normalization with
sequences of normalizing constants~$a(n) >0$ and $b(n)\in\mathbb{R}$, the sample maximum converges in distribution as:
\begin{equation}\label{ch1:maxdom}
\lim_{n\to\infty} P\Big(\frac{\max(X_1,\ldots,X_n)-b(n)}{a(n)}\leq x\Big) = \lim_{n\to\infty} F^n (a(n) x + b(n)) = G(x)\,,
\end{equation}
where $G$ is a non-degenerate distribution function which allows to gain limit theorems for random variables under consideration.

The class of distributions satisfying~\eqref{ch1:maxdom} is called the  \emph{max-domain of attraction} (MDA) of~$G$; when $F$ refers to this class we call it an \emph{MDA distribution}. Note that `essentially all the common continuous distributions of statistics or actuarial science are MDA distributions' (see \citealp{MFE05}, Sect.~7.1.2, p.~267).
A normalization of the maximum with $a(n)$ and $b(n)$ is necessary, because otherwise the limit distribution~$G$ would degenerate as $F^n (x)$ converges to zero for all $x$ smaller than the right endpoint $x_E := \sup\{x : F(x) <1\}$, or to one for all $x\geq x_E$ if $x_E$ is finite.
The prominent Fisher-Tippett theorem (see e.g. \citealp[Th.~3.2.3]{EKM97})
states that the limit distribution~$G$ belongs to the class of either Gumbel, Fr\'{e}chet, or Weibull distributions which are labelled as the \textit{extreme value distributions} characterized by the extreme value index~$\gamma$; more details are provided in Section~\ref{sec:2.1}.

However, already \citet[p.~67]{R87} noted that `computing normalizing constants $a(n), b(n)$ can be a brutal business'.
An elegant way to avoid computing normalizing sequences $a(n), b(n)$ in~\eqref{ch1:maxdom} is to consider -- instead of~\eqref{ch1:maxdom} -- alternative representations directly for MDA distributions~$F$ which are
characterized by a so-called \emph{auxiliary function}~$\psi$. In this paper we focus on this approach and investigate two popular representations of MDA distributions~$F$:
first, the variation representation (\VRp), where it holds that
\[\lim_{x\uparrow x_E} \frac{1-F(x+z \psi(x))}{1-F(x)} =
  \begin{cases}\exp(-z) \,, & $\!\!\!\!\!\!$ \text{if } \g= 0\\
(1+\gamma z)^{-1/\gamma}\,,\quad & $\!\!\!\!\!\!$ \text{if } \g\neq 0
  \end{cases}\]
for all $z$ with $1+\gamma z >0$, and second, the von Mises representation (\vMRp)
\[1-F(x) = c(x) \exp\Big(-\int_{x^\star}^x 1/\psi(u) \mathrm{d}u\Big)\,,\;\; x^\star < x < x_E\]
for some $x^\star$ and under certain regularity properties on functions~$c$ and~$\psi$; the exact definitions of \VR and \vMR are stated in Section~\ref{sec:2.2}.
\VR combines the concepts of both regular and gamma variation
and it is related to the peaks-over-threshold approach;
\vMR can be seen as a generalization of the important family of von Mises functions (cf. \citealp[Sect.~1.1]{R87});
for $\g\neq0$ the \vMR is also known as \emph{Karamata representation}.
Whether \VR or \vMR is to apply, is a matter of convenience: \vMR allows a direct relation between  distribution and auxiliary functions, whereas \VR provides an indirect relation between them which is defined via the ratio of survival functions.
These representations with corresponding auxiliary functions are also applicable in multivariate settings for modeling the asymptotic tail behavior of random vectors, e.g. for the Gumbel class with light or mildly heavy tails see \citet{H13}, \citet{NW21}, \citet{O17} or Seifert~(\citeyear{S14}, \citeyear{S16}).

Surprisingly, up to now there has been no literature with a systematic investigation concerning the validity of auxiliary functions for both representations for any of the MDA distribution classes.
In particular, the set of valid auxiliary functions for \vMR has neither been characterized completely nor separated
from those for \VRp, see e.g. \citet{EOW19}.
Moreover, in the current literature one could find statements such as `the auxiliary function for \VR is unique up to asymptotic equivalence' and `the auxiliary functions for both representations are asymptotically the same' (see e.g. Eq.~1.2.11 and Rem.~1.2.7 in \citealp{HF06}),
which could result in a misleading supposition that any auxiliary function which is valid for \VR could be also used for \vMRp.
In this paper we show that this is
in general not
the case: for example, \VR for distributions in the Gumbel max-domain of attraction is often applied with simple auxiliary functions, e.g. for Gaussian distribution with $\psi(x)=1/x$ or for Lognormal distribution with $\psi(x)=x/\log x$, however, as we show, these auxiliary functions are not valid for \vMR of these distributions. Furthermore, it is often insufficient only to know that an auxiliary function exists, but also the knowledge of the specific analytical form of $\psi$ is required. This could be the case when parameter estimation is of interest, see e.g. \citet[Eq.~14]{AFGS08}, \citet[Eq.~37]{FS12} who assume a certain form of~$\psi$ for estimation of the auxiliary function parameters.
In order to explore which analytical forms of~$\psi$ cover which subclasses of MDA distributions, it is desired to make a systematic analysis of \VR and \vMR with respect to their valid auxiliary functions~$\psi$.

In this paper we provide a comprehensive analysis of \VR and \vMR for the entire MDA distribution family.
Our major contributions are as follows: first, we prove that the set of \vMRp-valid auxiliary functions is a strict subset of
the set of \VRp-valid auxiliary functions, and show that some popular \VRp-valid functions are not \vMRp-valid. Second, we introduce a class of ``universal'' auxiliary functions $\psi_{\mathrm{u}}$ which are suitable for both \VR and \vMR for each of the three MDA classes.
Additionally, we provide explicit expressions for $\psi_{\mathrm{u}}$ which, however, require cumbersome integration. Next, we show under which assumptions $\psi_{\mathrm{u}}$ can be gained directly from the hazard functions. As expressions for $\psi_{\mathrm{u}}$  could still be analytically rather complicated, there arises a natural  question: how much can
one deviate from $\psi_{\mathrm{u}}$ to have still a valid auxiliary function $\psi$? To answer this question, we establish the sets of either \VRp- or \vMRp-valid auxiliary functions which provide ``if and only if'' characterizations, i.e. an auxiliary function is valid
``if and only if'' it refers to a corresponding set, respectively. This result is of much importance, because then for a given distribution function~$F$ it is sufficient to find a single valid universal auxiliary function in order to determine the whole set of valid auxiliary functions for this distribution. This characterization idea provides us a convenient toolkit for finding valid auxiliary functions of analytically simple form. Using this approach we deduce suitable (simple) auxiliary functions for several popular distributions for all three classes of the MDA distribution family.

Our paper is organized as follows. In Section~\ref{sec:2} we review the established results for MDA distributions and point on advantages of \VR and \vMR in contrast to applying directly the definition of MDA as in~\eqref{ch1:maxdom}. In Sections~\ref{sec:3} and~\ref{sec:4} we state our main novel results:
In Theorem~\ref{Th:MisesimpliesVar} we show that the set of auxiliary functions valid for \vMR is a strict subset of those valid for \VRp.
In Proposition~\ref{Prop:psi:lit} we prove that some popular \VRp-valid auxiliary functions remain \vMRp-valid only for subsets of MDA distributions. Next, in Theorem~\ref{Prop:zeta} we establish a class of universal auxiliary functions~$\psi_{\mathrm{u}}$ suitable for both representations.
Then in Proposition~\ref{Prop:psi0:2} we give under which additional assumptions the reciprocal hazard function is a valid choice for~$\psi_{\mathrm{u}}$; and in Proposition~\ref{Prop:psi0} we provide explicit expressions for $\psi_{\mathrm{u}}$ which could, however, be rather complicated. In order to find analytically more simple auxiliary functions,   in Theorem~\ref{Th:psi} we present a complete ``if and only if'' characterization of valid auxiliary functions for both representations of MDA distributions. In Section~\ref{sec:5} we discuss consequences of our results by providing suitable  auxiliary functions for some popular MDA distributions in Corollaries~\ref{Cor:EV} and~\ref{Cor:special-distr}. Moreover, in Section~\ref{sec:5} we show how our results can be applied for estimation of auxiliary function parameters.
Section~\ref{sec:6} concludes the paper, whereas the proofs are placed in Section~\ref{sec:7}.

\section{Variation and von Mises representations of MDA distributions}\label{sec:2}

\subsection{Preliminary results and notation}\label{sec:2.1}

Consider a random variable~$X$ with distribution function~$F$ in the max-domain of attraction (MDA) of some extreme value distribution~$G$, i.e. we study those~$F$ which satisfy~\eqref{ch1:maxdom}. Then due to the Fisher-Tippett theorem,
it follows that $G(x)=G_{\gamma}(p x + q)$ with some appropriate scale and location parameters $p>0$, $q\in\mathbb{R}$ so that
\begin{equation}\label{eq:EV}
G_{\gamma}(x)=\begin{cases} \exp(-\mathrm{e}^{-x})\,,\;\; x\in\mathbb{R}\,, & \text{if } \g= 0\\
\exp(-(1+\gamma x)^{-1/\gamma})\,,\;\; 1+\gamma x>0\,, \quad &\text{if } \g\neq 0
\end{cases}
\end{equation}
i.e. $G$~belongs to the type  $G_\g$ for some~$\gamma\in\mathbb{R}$. We write $F\in$~MDA$(G_{\gamma})$ with extreme value index~$\gamma$; then $G_{\gamma}$ is referred to as the \emph{generalized extreme value (GEV) distribution} with a single parameter~$\gamma$.
The family of GEV distributions~$G_\g$ in~\eqref{eq:EV} consists of three important classes: Gumbel for $\gamma=0$, Fr\'{e}chet for $\gamma>0$, and Weibull for $\gamma<0$:
\begin{eqnarray}\label{eq:EV:Def}
\begin{tabular}{llcll}
Gumbel-GEV for $\gamma=0$: $\!\!\!$ & $G_0(x)$ & $\!\!\!\!\!=\!\!\!\!\!$ &
    $\exp\left(-\text{e}^{-x}\right),\quad x\in\mathbb{R}$ & \\[4mm]
\vspace{2mm}
Fr\'{e}chet-GEV for $\gamma>0$: $\!\!\!$ & $G_{\g}(x)$ & $\!\!\!\!\!=\!\!\!\!\!$ &
    $\begin{cases} 0, & x\leq -1/\g\\ \exp(-(1+\gamma x)^{-1/\gamma}), & x>-1/\g \end{cases}$ & \\[8mm]
Weibull-GEV for $\gamma<0$: $\!\!\!$ & $G_{\g}(x)$ & $\!\!\!\!\!=\!\!\!\!\!$ &
    $\begin{cases} \exp(-(1-|\gamma| x)^{1/|\gamma|}), & x\leq 1/|\g|\\ 1, & x>1/|\g| \end{cases}\,.$ &\\[8mm]
\end{tabular}
\end{eqnarray}
Then $F\in$~MDA$(G_{\gamma})$ is said to be a Gumbel MDA distribution for $\gamma=0$, a Fr\'{e}chet MDA distribution for $\gamma>0$, or a Weibull MDA distribution for $\gamma<0$.
The Fr\'{e}chet MDA contains \textit{heavy-tailed} distributions, like Pareto, Cauchy, or Loggamma distributions. The Weibull MDA comprises \text{short-tailed} distributions with finite right endpoint $x_E <\infty$,
like Beta distributions.
The Gumbel MDA contains a large variety of distributions with rather different tail behavior: \textit{light-tailed} distributions like Gaussian, Gamma, Exponential, or light-tailed Weibull distributions (to be distinguished from the Weibull-GEV type) where the survival function~$\oF:=1-F$ decreases at least exponentially fast, as well as \textit{mildly heavy-tailed} distributions like Lognormal or heavy-tailed Weibull distributions where $\oF$ decreases faster than any power function but slower than an exponential function.
Note that these max-domains of attraction differ in the right endpoint $x_E = \sup\{x : F(x) <1\}\in\mathbb{R}\cup\{\infty\}$ of the distribution~$F$: for Fr\'{e}chet MDA it holds $x_E = \infty$, for Weibull MDA $x_E < \infty$, and for Gumbel MDA it can either be $x_E <\infty$ or $x_E=\infty$.

To describe the extreme tail behavior of distributions, it is not always advantageous to work directly with definition~\eqref{ch1:maxdom} because it might be a difficult task to find appropriate normalizing functions $a(n), b(n)$,
e.g. see \citet[p.67]{R87},
as it requires to determine the generalized inverse~$(1/\oF)^\leftarrow$ function.  To overcome the problem, \citet{R87} suggests to consider not~$F$ but a more simple tail-equivalent~$\tilde{F}$ and to compute~$a(n), b(n)$ for~$\tilde{F}$ instead. But even this way can lead to extensive calculations, see \citet[Example~2, pp.~71-72]{R87} for the Gaussian distribution.

When we confront a distribution referring to the Fr\'{e}chet class with $\gamma>0$, it is possible to apply a regular variation approach, see \citet[Sect.~0.4.1]{R87}. In this case, the survival function~$\oF$ can be represented as $\oF(x)=L(x) x^{-1/\g}$, $x>0$, with the variation index~$-1/\gamma$ and a slowly varying function~$L$ satisfying the condition \linebreak 
$\lim_{x\to\infty}L(\lambda x)/L(x)= 1$ for all $\lambda>0$. Such slowly varying functions are in general not easy to handle as they do not have to converge asymptotically, e.g. for Loggamma distribution which belongs to the Fr\'{e}chet MDA. A similar argumentation also applies to the Weibull class with $\gamma<0$.

For these reasons, instead of a slowly varying function~$L$ it appears to be advantageous to use a so-called auxiliary function $\psi$ which is often analytically more convenient compared to slowly varying functions. The use of auxiliary functions leads either to von Mises representation (\vMRp) or to variation representation (\VRp) which are in focus of this paper.

For description of the asymptotic behavior of functions we use the following notation: We write $f(x)\sim g(x)$, $x\to a$ for some $a \in \mathbb{R}\cup\{\infty\}$ in case when $f$ and $g$ are asymptotically equivalent functions, i.e. $\lim_{x\to a}f(x)/g(x)=1$.

\subsection{Representations for extreme value distributions and their MDA}\label{sec:2.2}

\noindent Now we define \VR and \vMR of the distribution function $F\in$~MDA$(G_\gamma)$ which are both based on the survival function $\oF=1-F$.
\begin{itemize}
\item[I.]
\textbf{Variation representation (\VRp):}
\begin{equation}\label{eq:repVar}
  F \in \text{MDA}(G_\gamma) \;\;\Leftrightarrow\;\; \lim_{x\uparrow x_E} \frac{\oF(x+z \psi(x))}{\oF(x)} =
  \begin{cases}\exp(-z) \,, & $\!\!\!\!\!\!$ \text{if } \g= 0\\
(1+\gamma z)^{-1/\gamma}\,,\quad & $\!\!\!\!\!\!$ \text{if } \g\neq 0
  \end{cases}
\end{equation}
for all $z$ satisfying $1+\gamma z >0$ with some eventually positive \emph{auxiliary function}~$\psi$;
thereby `eventually positive' means that there exists some $x^\star$ such that $\psi(x)>0$ for all $x>x^\star$.

\item[II.]
\textbf{Von Mises representation (\vMRp):}
\begin{equation}\label{eq:repMises}
  \!\! F \in \text{MDA}(G_\gamma) \;\;\Leftrightarrow\;\; \oF(x) = c(x) \exp\Big(-\!\int_{x^\star}^x 1/\psi(u) \mathrm{d}u\Big)\,,\; x^\star < x < x_E\,,
\end{equation}
for some $x^\star$ and function $c:(x^\star, x_E)\to(0,\infty)$ which converges to a positive limit, i.e. $\lim_{x\uparrow x_E} c(x)= c\in(0,\infty)$; and with auxiliary function $\psi$ satisfying \textbf{Property $\Pg$} which states:
\begin{eqnarray} \label{eq:Pg}
\hspace*{-7cm} \psi:(x^\star, x_E)\to(0,\infty) \text{ where}
\end{eqnarray}
$\quad$\\[-8mm]
\begin{tabular}{ll}
$\psi$ is absolutely continuous on $(x^\star,x_E)$ with $\lim_{x\uparrow x_E} \psi'(x)=0$\,, &$\!\!\!\!\!\!$ if $\;\gamma=0$\,,$\!\!$ \\
$\psi$ is continuous on $(x^\star,\infty)$ with $\lim_{x\to\infty} \psi(x)/x=\gamma$\,, &$\!\!\!\!\!\!$ if $\;\gamma>0$\,,$\!\!$ \\
$\psi$ is continuous  on $(x^\star,x_E)$ with $\lim_{x\uparrow x_E} \psi(x)/(x_E-x)=-\gamma$\,, &$\!\!\!\!\!\!$ if $\;\gamma<0$\,.$\!\!$\\
\end{tabular}
\end{itemize}

Both representations~\eqref{eq:repVar} and~\eqref{eq:repMises} are well-known and widely applied in the extreme value literature, see e.g. \citet[Th.~1.2.5, Th.~1.2.6]{HF06}.
We denote~\eqref{eq:repMises} as \emph{von Mises representation} (\vMRp) following \citet{J10} and \citet{S16} who interpreted it as a generalization of the representation of the so-called von Mises functions, cf. \citet[Eq.~1.3]{R87}; this representation was suggested initially for Gumbel MDA distributions (in a slightly different form) by \citet{BH72}.
In this paper we use the formulation of property~$\Pg$ from \citet[p.1354]{BH72} or \citet[Eqs.~1.3,1.5]{R87} for $\g=0$, and those from \citet[Th.~1.2.6]{HF06} for $\g\neq0$.
Further, we denote~\eqref{eq:repVar} as \emph{variation representation} (\VRp) as it covers both concepts of regular and gamma variation, cf. \citet[Sect. 0.4.1, 0.4.3]{R87}.

The variation representation is related to the peaks-over-threshold (POT) approach, an important technique in the statistics of extremal events. To illustrate this, note that for positive values~$z$ and sufficiently large~$x$ such that $\psi(x)$ is positive, the quotient $\oF(x+z \psi(x))/\oF(x)$ on the left hand side in the definition of \VR in~\eqref{eq:repVar} can be interpreted as excess probability
\begin{equation}
  P\left(\frac{X-x}{\psi(x)} > z\; \Big|\; X>x\right)\,,
\end{equation}
for a random variable~$X$ with the distribution function~$F$, and the limit on the right hand side in equation~\eqref{eq:repVar} is the so-called standard generalized Pareto distribution (GPD) for positive~$z$ with $1+\gamma z>0$.
In this sense, the statement in~\VR gives
a distributional approximation for the scaled excesses over a high threshold~$x$.
This is related to the idea behind the popular POT approach, which grounds on the approximation of the excess distribution function by the GPD; more details can be found e.g. in \citet[Sect.~3.4, 6.5]{EKM97} and \citet{P75}.

\VR and \vMR are rather advantageous for gaining limit results.
Both~\eqref{eq:repVar} and~\eqref{eq:repMises} help to get statements about the asymptotic tail behavior of~$\oF$ by a \emph{single} auxiliary function~$\psi$ which could be determined conveniently for many popular distributions. This is because \VR and \vMR allow for a certain degree of freedom in the choice of auxiliary functions~$\psi$, as we show in Theorem~\ref{Th:psi} of this paper. Also the tail estimation often grounds on these representations based on auxiliary function~$\psi$, see e.g. \citet{BFR17}, \citet{AFGS08}. Note that \vMR provides a direct expression of~$\oF$ which is sometimes advantageous compared to indirect definitions of $\oF$ in~\eqref{ch1:maxdom}, or in~\eqref{eq:repVar} for \VRp.

\VR and \vMR apply for all MDA distributions with arbitrary $\g\in\mathbb{R}$, not only for $\g\neq 0$\linebreak
 as it is the case in regular variation approach.
In particular, \VR and \vMR are often used for Gumbel MDA distributions, see \citet{BS16}, \citet{DFH18}, \citet{FS12}, \citet{H12}, \citet{KL05} and \citet{KS20} for \VRp, as well as \citet{AFGS08}, \citet{EOW19}, \citet{FKL06}, \citet{J10} and \citet{S16} for \vMRp; this Gumbel case grants the most interesting analysis and results, as we will point out in the following sections. Using \VR and \vMR could be also convenient instead of applying the concept of regular variation suitable for Fr\'{e}chet and Weibull MDA distributions characterized by parameter $\gamma\neq 0$ (cf. e.g. \citealp{R87}; \citealp[Th 1.2.1(i)]{HF06}).
For example, we show that for \VR of Fr\'{e}chet or Weibull MDA distributions we can simply choose $\psi(x)=\gamma x$ or $\psi(x)=-\g (x_E -x)$, respectively, which are much more simple functions than those slowly varying functions~$L$ for regular variation, which is discussed in Section~\ref{sec:2.1}.

To sum it up, \VR and \vMR in~\eqref{eq:repVar} and~\eqref{eq:repMises}, respectively, are rather practical compared to alternative ways to represent MDA distributions. Consequently, it is of interest to provide a deeper insight concerning suitable auxiliary functions~$\psi$ and their properties and to gain a complete characterization for the class of valid auxiliary functions for both \VR and \vMRp, which is done in this paper.
The knowledge of the functional form of~$\psi$ is also very convenient for estimation purposes in application of extreme value results. Then $\psi$ can be represented as an analytical function with certain unknown parameters and estimation of $\psi$ reduces to estimation of these parameters, see e.g. \citet{AFGS08}, \citet{FS12}.
In Section~\ref{sec:5} we show how our results help to find auxiliary functions suitable for particular MDA distributions, which is important for estimation purposes.

\section{Universal auxiliary functions for \VR and \vMRs}\label{sec:3}

Now we provide our theoretical findings on valid auxiliary functions $\psi$ for distributions in the max-domain of attraction for all three MDA classes. More precisely, we present in Theorem~\ref{Th:MisesimpliesVar} the relation between auxiliary functions valid for \VR or for \vMRp, whereas in Theorem~\ref{Prop:zeta} we provide universal auxiliary functions valid for both \VR and \vMRp, which is one of our central results.

\begin{theorem}\label{Th:MisesimpliesVar}
Let $F$ be a distribution function in the max-domain of attraction MDA($G_{\gamma}$) of some GEV distribution $G_\g$ with arbitrary $\gamma\in\mathbb{R}$.
If $\psi$ is an auxiliary function valid for \vMR in~\eqref{eq:repMises}, then $\psi$ is also valid for \VR in~\eqref{eq:repVar}.
\end{theorem}

The result in Theorem~\ref{Th:MisesimpliesVar} implies that in order to find universal auxiliary functions suitable for both \VR and \vMR it is sufficient to analyze only \vMR with respect to valid auxiliary functions.
However, it is not sufficient for a \VRp-valid function to satisfy only the property~$\Pg$ in~\eqref{eq:Pg} in order to be a \vMRp-valid auxiliary function.
As a counter-example, in case of Loggamma distribution with positive parameters $\alpha$ and $\beta$, the function $\psi(x)=x/\alpha$ obviously satisfies property~$\Pg$, $\g>0$, and is also valid for \VRp, but (except for the special case of a Pareto distribution with $\beta=1$) this $\psi$ is \emph{not} valid for \vMRp,
which is proven in Remark~\ref{Rem:Examples:Sect3}(ii).

Next we re-consider some popular \VRp-valid auxiliary functions frequently mentioned in the extreme value literature (cf. \citealp[Th. 1.2.1, 1.2.5]{HF06}; \citealp[p.143, Eq.~3.34]{EKM97}; \citealp[Prop. 1.9]{R87})
for each of the three MDA classes:
\begin{itemize}
\item for Gumbel MDA with $\gamma=0$:
\begin{eqnarray}\label{eq:psi:Gum:1}
\psi(x) &=& \frac{\int_x^{x_E}\oF(u)\mathrm{d}u}{\oF(x)}\,,\\ \label{eq:psi:Gum:2}
\psi(x) &=& \frac{\int_x^{x_E}\int_v^{x_E}\oF(u)\mathrm{d}u\mathrm{d}v}{\int_x^{x_E}\oF(u)\mathrm{d}u}\,,
\end{eqnarray}
\item for Fr\'{e}chet MDA with $\gamma>0$:
\begin{equation}\label{eq:psi:Fre}
\psi(x) = \gamma x\,,
\end{equation}
\item for Weibull MDA with $\gamma<0$:
\begin{equation}\label{eq:psi:Wei}
\psi(x) = -\gamma (x_E -x)\,.
\end{equation}
\end{itemize}

However, for most distributions $F$ in MDA$(G_\g)$ the functions $\psi$ from \eqref{eq:psi:Gum:1}--\eqref{eq:psi:Wei} do not define valid auxiliary functions for \vMRp, which is illustrated by the examples of the Gaussian distribution (as a representant of the Gumbel MDA class) and the Loggamma distribution (as a representant of the Fr\'{e}chet MDA class) in the following Remark:

\begin{remark}\label{Rem:Examples:Sect3} $\quad$\\
(i)$\;\;$ First we show that choice $\psi$ from~\eqref{eq:psi:Gum:1} is not valid for standard Gaussian distribution $F=\Phi$. For $x>x^\star$ with arbitrary $x^\star$ it holds that:
\begin{eqnarray*}
&&\exp\Big(-\int_{x^\star}^x \frac{\overline{\Phi}(v)}{\int_v^{\infty}\overline{\Phi}(u)\mathrm{d}u} \mathrm{d}v\Big)
=  \exp\Big(-\int_{x^\star}^x \Big(-\frac{\mathrm{d}}{\mathrm{d}v}\log \int_v^{\infty}\overline{\Phi}(u) \mathrm{d}u\Big) \mathrm{d}v\Big)\\
&=& \exp\Big(-\Big(-\log \int_x^{\infty}\overline{\Phi}(u) \mathrm{d}u + \log \int_{x^\star}^{\infty}\overline{\Phi}(u) \mathrm{d}u\Big)\Big) = \frac{\int_x^{\infty}\overline{\Phi}(u) \mathrm{d}u}{k({x^\star})}\,,
\end{eqnarray*}
with $k({x^\star}):=\int_{x^\star}^{\infty}\overline{\Phi}(u) \mathrm{d}u \in(0,\infty)$.
If $\psi$ from~\eqref{eq:psi:Gum:1} would be a \vMRp-valid auxiliary function for the Gaussian distribution, then $\overline{\Phi}$ would have to satisfy \vMR in~\eqref{eq:repMises} where $c(x)=k(x^\star)\, \overline{\Phi}(x)/\int_x^{\infty} \overline{\Phi}(u) \mathrm{d}u$ would converge to a positive finite limit for $x\to\infty$. But the Mill's ratio $\overline{\Phi}(x)\sim \varphi(x)/x$ with standard normal density~$\varphi$ implies that
\[\lim_{x\to\infty} \frac{\int_x^{\infty} \overline{\Phi}(u) \mathrm{d}u}{\overline{\Phi}(x)} = \lim_{x\to\infty} \frac{\overline{\Phi}(x)}{\varphi(x)}=0\,.\]
Hence, function $c$ does not possess a finite limit and, consequently, $\psi$ from~\eqref{eq:psi:Gum:1} is \emph{not} a \vMRp-valid auxiliary function for Gaussian distribution.\\
(ii)$\;\;$ Next, we show that choice $\psi=\g x$ from~\eqref{eq:psi:Fre} is not \vMRp-valid for Loggamma distribution with $\oF(x) = \left(\alpha^{\beta}/\Gamma(\beta)\right)\, \int_{x}^{\infty} \left( u^{-\alpha-1} (\log u)^{\beta-1}\right)\mathrm{d}u$, $x>1$, $\alpha, \beta>0$, where $\beta\neq1$, which is a Fr\'{e}chet MDA distribution with extreme value index~$\gamma=1/\alpha$:
For $x>x^\star$ it holds that $\exp(-\int_{x^\star}^x \alpha/u\, \mathrm{d}u) = (x^\star)^{\alpha}\, x^{-\alpha}$.\\
If $\psi$ from~\eqref{eq:psi:Fre} would be a \vMRp-valid auxiliary function for Loggamma distribution, then $\oF$ would have to satisfy \vMR in~\eqref{eq:repMises} where $c(x)=(x^\star)^{-\alpha} x^{\alpha} \oF(x)$ would converge to a positive finite limit for $x\to\infty$. But for Loggamma it holds that
\[\oF(x)\sim \frac{\alpha^{\beta-1}}{\Gamma(\beta)} x^{-\alpha} (\log x)^{\beta-1}\,,\; x\to\infty\,.\]
Hence, for $\beta\neq1$, function $c$ does not possess a limit in $(0,\infty)$ and, consequently, $\psi$ from~\eqref{eq:psi:Fre} is \emph{not} a \vMRp-valid auxiliary function for Loggamma distribution, except for the special case of a Pareto distribution for $\beta=1$.
$\hfill \diamond$
\end{remark}

Now we formulate conditions on distributions~$F$ required under which the \VRp-valid functions~$\psi$ from \eqref{eq:psi:Gum:1}--\eqref{eq:psi:Wei} are also \vMRp-valid.

\begin{proposition}\label{Prop:psi:lit}
Let $F$ be a distribution function in the max-domain of attraction MDA($G_{\gamma}$) of some GEV distribution $G_\g$ with arbitrary $\gamma\in\mathbb{R}$.
Then it holds:
\begin{itemize}
\item for Gumbel MDA with $\gamma=0$: $\;\psi$ in~\eqref{eq:psi:Gum:1} is a valid auxiliary function for \vMR if and only if $F$ is (eventually) absolutely continuous and has an exponential tail behavior with infinite right endpoint $x_E=\infty$, i.e. $\oF(x)\sim K\, \exp(-\lambda x + h(x))$ as $x\to\infty$ with $K, \lambda>0$ and twice differentiable~$h$ with $\lim_{x\to\infty}h'(x)=\lim_{x\to\infty} h''(x)=0$. The result for $\psi$ in~\eqref{eq:psi:Gum:2} is similar, however, the assumption about the absolute continuity of~$F$ is not required in that case;
\item for Fr\'{e}chet MDA with $\gamma>0$: $\;\psi$ in~\eqref{eq:psi:Fre} is a valid auxiliary function for \vMR if and only if $F$ follows a Pareto-like distribution, i.e. $\oF(x)\sim K\, x^{-1/\g}$ with $K>0$ and $x\to\infty$;
\item for Weibull MDA with $\gamma<0$: $\;\psi$ in~\eqref{eq:psi:Wei} is a valid auxiliary function for \vMR if and only if $F$ follows a Pareto-like distribution with finite right endpoint $x_E<\infty$, i.e. $\oF(x)\sim K\, (x_E-x)^{-1/\g}$ with $K>0$ and $x\uparrow x_E$.
    \end{itemize}
\end{proposition}

Proposition~\ref{Prop:psi:lit} states that functions~\eqref{eq:psi:Fre} or~\eqref{eq:psi:Wei} define \vMRp-valid auxiliary functions only for Pareto-like distributions, but not for other important Fr\'{e}chet or Weibull MDA distributions.
Further, some Gumbel MDA members like exponential or Gamma distributions have exponential tail behavior as required in Proposition~\ref{Prop:psi:lit} for $\g=0$, but many important Gumbel MDA distributions (e.g. Gaussian, Lognormal and generic Weibull distributions) do not have exponential tail behavior and, hence, functions  in~\eqref{eq:psi:Gum:1} and~\eqref{eq:psi:Gum:2} are not \vMRp-valid auxiliary functions.
To sum it up, none of the \VRp-valid auxiliary functions in \eqref{eq:psi:Gum:1}--\eqref{eq:psi:Wei} is in general valid for \vMRp.

An alternative approach to find valid auxiliary functions for Gumbel MDA distributions is proposed by \citet{HF06} who state in their Remark~1.2.7 that for Gumbel MDA `the auxiliary functions [for both representations] are asymptotically the same' and suggest to choose
\begin{equation} \label{eq:psi:Gum:3}
\psi(x) = \frac{\oF(x)}{F'(x)}
\end{equation}
when it holds that (i) $F$ is a Gumbel MDA distribution; (ii) $F$ is twice differentiable with a positive first derivative~$F'$; (iii) $F$ satisfies the so-called von Mises condition:
\begin{equation} \label{eq:psi:vMcond}
\lim_{x\uparrow x_E}\frac{\oF(x) F''(x)}{(F'(x))^2}=-1\,.
\end{equation}
However, the condition~\eqref{eq:psi:vMcond} is not satisfied for all Gumbel MDA distributions.
This approach has already been stated in the early literature, e.g. \citet[Prop.~1.1(b)]{R87}, \citet[Ex.~3.3.23]{EKM97}.

Altogether, up to now in the current literature there has been no discussion about how to choose \vMRp-valid auxiliary functions. Here we close this gap by
establishing a form of~$\psi$, such that it is \vMRp-valid
for all distributions referring to either Gumbel, Fr\'{e}chet, or Weibull MDA classes.
Our approach grounds on the following idea:
for two continuous functions $\zeta$ and $\psi_{\mathrm{u}}$ whose connection is established by
\begin{equation}\label{eq:zeta-psi0}
\frac{\mathrm{d}}{\mathrm{d}x} \log \zeta(x) = -\frac{1}{\psi_{\mathrm{u}}(x)}\,,
\end{equation}
it holds for $x>z$ that:
\begin{equation}\label{eq:zeta-psi0:1}
\zeta(x) = c_z \exp\Big(-\int_z^x 1/\psi_{\mathrm{u}}(u)\mathrm{d}u\Big)\,,
\end{equation}
for some $z$ and a positive constant $c_z$. The last expression is obtained from equation~\eqref{eq:zeta-psi0} by integrating both sides and solving it for~$\zeta$.
Hence, the issue of searching for a \vMRp-valid auxiliary function can be reduced to finding an absolutely continuous asymptotic equivalent function~$\zeta$ for $\oF$, i.e. $\zeta(x)\sim \oF(x)$ as $x\uparrow x_E$. We obtain from equation~\eqref{eq:zeta-psi0} that when
\begin{equation}\label{eq:zeta-psi0:2}
  \psi_{\mathrm{u}}(x)=-\zeta(x)/\zeta'(x)
\end{equation}
satisfies the property~$\Pg$, i.e. $ \psi_{\mathrm{u}}\in\Pg$ with $\Pg$ from~\eqref{eq:Pg}, then $\psi_{\mathrm{u}}$ is a valid universal auxiliary function for both \vMR and \VRsp. This leads us to the following result:

\begin{theorem}\label{Prop:zeta}
  Let $F$ be a distribution function in the max-domain of attraction MDA($G_{\gamma}$) of some GEV distribution $G_\g$ with arbitrary $\gamma\in\mathbb{R}$.
  Then there exists an absolutely continuous function~$\zeta$ with $\oF(x)\sim \zeta(x)$ as $x\uparrow x_E$ which satisfies that $-\zeta/\zeta'\in\Pg$ as in~\eqref{eq:Pg}, and the function
  \begin{equation}\label{eq:zeta:Th}
    \psi_{\mathrm{u}} (x):=-\frac{\zeta(x)}{\zeta'(x)}
  \end{equation}
  is a universal auxiliary function valid for both \VR in~\eqref{eq:repVar} and \vMR in~\eqref{eq:repMises}.
\end{theorem}

The result in Theorem~\ref{Prop:zeta} offers a variety of possibilities to obtain universal auxiliary functions~$\psi_{\mathrm{u}}$ which are valid for both \VR and \vMRp. Such an approach is rather advantageous because it allows:
\begin{itemize}
\item[(A)] to prove in Proposition~\ref{Prop:psi0:2} that under some additional assumptions on distribution~$F$ referring to any of the three MDA classes, the reciprocal hazard function $\psi_{\mathrm{u}}=\oF/F'$  defines a valid universal auxiliary function for $F$;
\item[(B)] to obtain in Proposition~\ref{Prop:psi0} explicit expressions for universal auxiliary functions for any distribution referring to  each of the three MDA classes;
\item[(C)] to derive convenient universal auxiliary functions for several important distributions, which are analytically more simple than those based on the reciprocal hazard function $\oF/F'$; we  provide examples of such convenient auxiliary functions in the table referring to Corollary~\ref{Cor:special-distr}.
\end{itemize}

To address (A), in the following Proposition~\ref{Prop:psi0:2} we provide additional assumptions to be satisfied by~$F$ which guarantee that the immediate candidate~$\zeta=\oF$ satisfies the requirements of Theorem~\ref{Prop:zeta}.
This proposition extends the result of \citet{HF06} for Gumbel MDA, cf.~\eqref{eq:psi:Gum:3}, on all three MDA classes.

\begin{proposition}\label{Prop:psi0:2}
Let $F$ be a distribution function in the max-domain of attraction MDA($G_{\gamma}$) of some GEV distribution $G_\g$ with arbitrary $\gamma\in\mathbb{R}$.
Let $F$ be absolutely continuous in some left neighborhood of $x_E$ with positive and non-increasing density~$F'$. If $\gamma=0$, let additionally assume that $F$ is twice differentiable in some left neighborhood of $x_E$ with negative and non-decreasing~$F''$, and let $F'(x)=\int_x^{x_E}(-F''(u))\mathrm{d}u$ for $x$ close to $x_E$.
Under these assumptions, it follows that the reciprocal hazard function, namely
\begin{equation}\label{eq:prop:psi0:add}
\psi_{\mathrm{u}}(x)=\frac{\oF(x)}{F'(x)}\,,
\end{equation}
is a valid auxiliary function for both \VR in~\eqref{eq:repVar} and \vMR in~\eqref{eq:repMises} of the survival function~$\oF$.
Moreover, these assumptions imply that $F$ has a \vMR with constant function $c(x)\equiv c\in(0,\infty)$.
\end{proposition}

Note that the \vMR with constant $c(x)\equiv c$ do not cover the entire MDA distribution family.
For this reason, as noted in~(B), in Proposition~\ref{Prop:psi0} we apply our result from Theorem~\ref{Prop:zeta} and find functions~$\zeta$ which provide us universal auxiliary functions $\psi_{\mathrm{u}}$ valid for \emph{all}  distributions in each of the three MDA classes.

\begin{proposition}\label{Prop:psi0}
Let $F$ be a distribution function in the max-domain of attraction MDA($G_{\gamma}$) of some GEV distribution $G_\g$ with arbitrary $\gamma\in\mathbb{R}$. Depending on the MDA class, we define universal auxiliary functions $\psi_{\mathrm{u}}$ as follows:\\
$\bullet\;$ for Gumbel MDA with $\gamma=0$:
\begin{equation}\label{eq:psi0:G}
\psi_{\mathrm{u}}(x) := \frac{\II \III }{3 \I \III -2 \left(\II\right)^2}\,,
\end{equation}
$\bullet\;$ for Fr\'{e}chet MDA with $\gamma>0$:
\begin{equation}\label{eq:psi0:F}
\psi_{\mathrm{u}}(x) := \frac{x\,\int_x^{\infty}(\oF(u)/u)\; \mathrm{d}u}{\oF(x)}\,,
\end{equation}
$\bullet\;$ for Weibull MDA with $\gamma<0$:
\begin{equation}\label{eq:psi0:W}
\psi_{\mathrm{u}}(x) := \frac{(x_E-x)\int_{x}^{x_E}\oF(u)/(x_E-u)\,\mathrm{d}u}{\oF(x)}\,.
\end{equation}
Then it holds that these functions $\psi_{\mathrm{u}}$ are valid auxiliary functions for both  \VR in~\eqref{eq:repVar} and \vMR in~\eqref{eq:repMises} of the survival function~$\oF$.
\end{proposition}

In Proposition~\ref{Prop:psi0} we provide the explicit expressions for universal auxiliary function suitable for all distributions referring to one of the three MDA classes.
However, the price to pay for the generality of the result in Proposition~\ref{Prop:psi0} is a possibly circumstantial analytical form of these auxiliary functions, in particular for the Gumbel MDA class. Further, it is also not always easy to exploit the results of Propositions~\ref{Prop:psi0:2} and~\ref{Prop:psi0}, because the explicit analytical form of the survival function $\oF$ is not available for many important distributions. Fortunately, our Theorem~\ref{Prop:zeta} provides a method to construct (by using appropriate functions~$\zeta$) more simple auxiliary functions even for distributions satisfying~\eqref{eq:psi:vMcond}. For these reasons, in the following sections we explore (C) and show the way how to overcome these difficulties.
Moreover, in the following sections we show how much more freedom we have in the choice of auxiliary functions valid only for \VR and how we can simplify the analytical form of auxiliary functions if one is interested only in \VR and not in \vMR representations.

\section{Sets of \VRp- and \vMRp-valid auxiliary functions}\label{sec:4}

In the previous section we have derived the requirements for universal, i.e. both \VRp- and \vMRp-valid, auxiliary functions $\psi_{\mathrm{u}}$.
Although we can find such functions valid for all distributions in each of MDA classes by applying the results in Proposition~\ref{Prop:psi0}, they are not always easily analytically tractable. In order to find more simple auxiliary functions, we need to answer the question how much we can deviate from~$\psi_{\mathrm{u}}$ to have still a valid auxiliary function?
In the following theorem we state one of our major results which defines the complete sets of valid auxiliary functions for \VR or for \vMR such that we obtain an ``if and only if'' characterization. This means, an auxiliary function is either \VRp- or \vMRp-valid ``if  and only if'' it refers to the corresponding set, respectively.

\begin{theorem}\label{Th:psi}
Let $F$ be a distribution function in the max-domain of attraction MDA($G_{\gamma}$) of some GEV distribution $G_\g$ with arbitrary $\gamma\in\mathbb{R}$.
Given a universal auxiliary function $\psi_{\mathrm{u}}$ constructed according to Theorem~\ref{Prop:zeta} (e.g. the explicit $\psi_{\mathrm{u}}\,$ from \eqref{eq:prop:psi0:add}--\eqref{eq:psi0:W} as in Propositions~\ref{Prop:psi0:2} and~\ref{Prop:psi0}), it holds that:
\begin{itemize}
  \item[(i)] $\psi$ is a valid auxiliary function for \VR in~\eqref{eq:repVar} of~$F$ if and only if $\psi\in\AVF$ where
      \begin{equation}\label{eq:var:equiv}
   \AVF := \{\psi \,:\, \psi(x) \sim \psi_{\mathrm{u}}(x)\;\;\text{ as } x\uparrow x_E\}\,,
   \end{equation}
      which is the entire asymptotic equivalence set of the universal auxiliary function~$\psi_{\mathrm{u}}$.
  \item[(ii)] $\psi$ is a valid auxiliary function for \vMR in~\eqref{eq:repMises} of~$F$ if and only if $\psi\in\AMF$ where
    \begin{eqnarray} \nonumber
   \AMF &:=& \Big\{\psi \in\Pg\,:\, \exists\, z<x_E, K_z\in(-\infty,\infty) \Big. \qquad\qquad \\ \label{eq:int:h}
   && \Big. \;\;\;\qquad \text{\, with \,} K_z =\int_z^{x_E}\frac{\psi_{\mathrm{u}}(v) - \psi(v)}{\psi_{\mathrm{u}}(v)\psi(v)}\mathrm{d}v\,\Big\}\,,
   \end{eqnarray}
   with the property $\Pg$ defined in~\eqref{eq:Pg}.
\end{itemize}
Further, it holds that $\AMF$ is a strict subset of $\AVF$, i.e. $\AMF\subsetneq \AVF$.
\end{theorem}

The result in Theorem~\ref{Th:psi} is of much importance: it states that for a given distribution~$F$ it is sufficient to find a single universal auxiliary function~$\psi_{\mathrm{u}}$ in order to be able to specify the complete (``if and only if'') sets of valid auxiliary functions -- $\AVF$ or $\AMF$ -- for this distribution. Alternatively, for a given candidate function~$\psi$ it could be checked whether it refers to the sets of valid auxiliary functions for a given distribution.
Theorem~\ref{Th:psi} characterizes the freedom we have in the choice of auxiliary functions by specifying how much we can deviate from universal~$\psi_{\mathrm{u}}$ to get still a valid~$\psi$: for \VR the entire asymptotic equivalence set $\{\psi\sim\psi_{\mathrm{u}}\}$ is possible, whereas for \vMR the choice of valid~$\psi$ is much more restrictive as specified in~\eqref{eq:int:h}.
As any universal
auxiliary function $\psi_{\mathrm{u}}$ chosen according to Theorem~\ref{Prop:zeta} is suitable for these purposes,
we are able to find rather simple auxiliary functions valid for several important distributions in Corollaries~\ref{Cor:EV} and~\ref{Cor:special-distr}.

The representations of the classes $\AVF$ and $\AMF$ in Theorem~\ref{Th:psi} support the statement of Theorem~\ref{Th:MisesimpliesVar} that all \vMRp-valid functions~$\psi$ are also \VRp-valid.
Moreover, in Theorem~\ref{Th:psi} we prove that the opposite statement is in general not true: there exist functions~$\psi\sim\psi_{\mathrm{u}}$ which do not satisfy the finiteness condition of the integral in $\AMF$,
their examples are provided in Corollaries~\ref{Cor:EV} and~\ref{Cor:special-distr}. Furthermore, an explicit application of Theorem~\ref{Th:psi}(ii) as a criterion to check whether a candidate $\psi$ is \vMRp-valid is given in Remark~\ref{Rem:AMF:Criterion} in Section~\ref{sec:5}.

The universal auxiliary functions $\psi_{\mathrm{u}}$ given in Eqs.~\eqref{eq:psi0:G}--\eqref{eq:psi0:W} are defined in a quite different way for the Gumbel ($\gamma=0$), Fr\'{e}chet ($\gamma>0$), or Weibull ($\gamma<0$) MDA classes; the same is valid for the \VRp-valid auxiliary functions given in Eqs.~\eqref{eq:psi:Gum:1}--\eqref{eq:psi:Wei}. Remarkably, for all auxiliary functions~$\psi$ holds that the ratio $\psi(x)/x$ (or $\psi(x)/(x-x_E)$ for finite~$x_E$, respectively) converges to limit~$\g$ as $x\uparrow x_E$, as we show in the next corollary.

\begin{corollary}\label{Cor:psi0}
All auxiliary functions~$\psi$ valid for \VRp, i.e. all  $\psi\in\AVF$ and, consequently, all $\psi\in\AMF$ for arbitrary MDA distribution~$F$ satisfy the following property:
\[\lim_{x\to \infty} \frac{\psi(x)}{x} = \gamma\;\; \text{ if } x_E=\infty \quad \text{ or } \quad \lim_{x\uparrow x_E} \frac{\psi(x)}{x-x_E} = \gamma \;\; \text{ if } x_E<\infty\]
for arbitrary extreme value index $\gamma\in\mathbb{R}$.
\end{corollary}

\section{Consequences of the results}\label{sec:5}

Now we exploit our results derived in previous sections in order to make some recommendations concerning the choice of convenient valid auxiliary functions. In particular, we provide auxiliary functions suitable for several popular distributions, whereby we deduce universal auxiliary functions~$\psi_{\mathrm{u}}$ valid for both representations as well as parsimonious auxiliary functions~$\psi$ valid for \VR only.

First we apply the result in Theorem~\ref{Prop:zeta} to obtain convenient auxiliary functions for some important distributions. We start with the GEV distributions $G_\g$ from~\eqref{eq:EV}. In this case we could apply the result in Proposition~\ref{Prop:psi0:2} and use $\psi_{\mathrm{u}}=\oG_\g/G'_\g$ as it satisfies the property~$\Pg$, however, Theorem~\ref{Prop:zeta} allows us to find even more simple auxiliary functions by choosing appropriate asymptotically equivalent functions~$\zeta$:

\begin{itemize}
\item for Gumbel-GEV~$G_{0}$ it holds that:
\begin{equation}
  \oG_{0}(x) = 1-\exp(-\mathrm{e}^{-x}) \sim \mathrm{e}^{-x} =: \zeta(x)\,,\;\; x\to\infty\,,
\end{equation}
using that $1-\exp(-y)\sim y$ for $y\to0$.
\item for Fr\'{e}chet-GEV~$G_{\gamma}$, $\gamma>0$, it holds that:
\begin{equation}
  \oG_{\gamma}(x) = 1-\exp(-(1+\gamma x)^{-1/\gamma}) \sim (\gamma x)^{-1/\gamma} =: \zeta(x)\,,\;\; x\to\infty\,,
\end{equation}
whereby we use a power series expansion of the exponential function.
\item for Weibull-GEV~$G_{\gamma}$, $\gamma<0$ and right endpoint $x_E=-1/\gamma$, it holds that:
\begin{equation}
  \oG_{\gamma}(x) = 1-\exp(-(1+\gamma x)^{-1/\gamma}) \sim (-\gamma (x_E-x))^{-1/\gamma} =: \zeta(x)\,,\;\; x\uparrow x_E\,.
\end{equation}
\end{itemize}

As all these choices of $\zeta$ satisfy $\psi_{\mathrm{u}}(x):=-\zeta(x)/\zeta'(x)\in\Pg$, we now suggest simple universal auxiliary functions for the GEV in the following corollary:

\begin{corollary}\label{Cor:EV}
For the generalized extreme value distributions~$G_{\gamma}$ we define universal auxiliary functions
$\psi_{\mathrm{u}}$
depending on the extreme value index~$\gamma\in\mathbb{R}$ as follows:\\[1mm]
\begin{tabular}{llcll}
Gumbel-GEV  $G_0$: \hspace{10mm} & $\psi_{\mathrm{u}}(x)$ & $\equiv$ & $1$, & \\
Fr\'{e}chet-GEV $G_\gamma$, $\gamma>0$: & $\psi_{\mathrm{u}}(x)$ & $=$ & $\gamma x$, & \\
Weibull-GEV  $G_\gamma$, $\gamma<0$: & $\psi_{\mathrm{u}}(x)$ & $=$ & $-\gamma (x_E-x)$. &\\[1mm]
\end{tabular}

\noindent All these $\psi_{\mathrm{u}}$ are valid universal auxiliary functions for both \VR in~\eqref{eq:repVar} and \vMR in~\eqref{eq:repMises} for the GEV distributions from~\eqref{eq:EV}.
\end{corollary}

Next, we give some convenient auxiliary functions for important MDA distributions by finding suitable asymptotically equivalent functions~$\zeta$ based on the idea of~Theorem~\ref{Prop:zeta}. The choice of these distributions is motivated by results in financial literature, see e.g. \citet{EKM97}, \citet{GOS12}, \citet{KS19}.

\begin{itemize}
\item[] \hspace{-7mm} Gumbel MDA distributions:\\
  $\bullet\;$ (Standard) Gaussian distribution:
\begin{equation}\label{eq:specdist:A}
  \overline{\Phi}(x) \sim \frac{\varphi(x)}{x} = \frac{1}{\sqrt{2\pi} x} \exp(-x^2/2) =: \zeta(x)\,,\;\; x\to\infty\,;
\end{equation}
  $\bullet\;$ (Standard) Lognormal distribution
\begin{eqnarray} \nonumber
  \oF(x) &\!\!\!\!\!=\!\!\!\!\!& \overline{\Phi}(\log x)\\
  &\!\!\!\!\!\sim\!\!\!\!\!& \frac{\varphi(\log x)}{\log x} = \frac{1}{\sqrt{2\pi} \log x} \exp(-(\log x)^2/2) =: \zeta(x)\,,\;\; x\to\infty\,;
\end{eqnarray}
  $\bullet\;$ Exponential and Exponential-like distribution with $K, \lambda>0$:
\begin{equation}
  \oF(x) \sim K\, \exp(-\lambda x) =:\zeta(x)\,,\;\; x\to\infty\,;
\end{equation}
  $\bullet\;$ Gamma distribution with $x>0$ and $\alpha, \beta>0$:
\begin{eqnarray} \nonumber
  \oF(x) &=& \frac{\beta^{\alpha}}{\Gamma(\alpha)}\, \int_{x}^{\infty} \left( u^{\alpha-1} \exp(-\beta u)\right)\mathrm{d}u\\
   &\sim& \frac{1}{\Gamma(\alpha)}\, \left( \beta x\right)^{\alpha-1} \exp(-\beta x) =: \zeta(x) \,,\;\; x\to\infty\,;
\end{eqnarray}
  $\bullet\;$ Weibull (with $\alpha=0$) and Weibull-like distribution with $K, \beta, \tau>0$, $\alpha\in\mathbb{R}$:
\begin{equation}
  \oF(x) \sim K\, x^{\alpha}\, \exp(- (\beta x)^{\tau}) =: \zeta(x) \,,\;\; x\to\infty\,;
\end{equation}

\item[] \hspace{-7mm} Fr\'{e}chet MDA distributions:\\
  $\bullet\;$ Loggamma distribution with $x>1$ and $\alpha, \beta>0$:
\begin{eqnarray} \nonumber
  \oF(x) &=& \frac{\alpha^{\beta}}{\Gamma(\beta)}\, \int_{x}^{\infty} \left( u^{-\alpha-1} (\log u)^{\beta-1}\right)\mathrm{d}u\\
  &\sim& \frac{\alpha^{\beta-1}}{\Gamma(\beta)} x^{-\alpha} (\log x)^{\beta-1} =: \zeta(x) \,,\;\; x\to\infty\,;
\end{eqnarray}
  $\bullet\;$ (Standard) Cauchy distribution:
\begin{equation}
  \oF(x) = \int_x^{\infty} \frac{1}{\pi (1+u^2)} \mathrm{d}u \sim \frac{1}{\pi x} =:\zeta(x)\,,\;\; x\to\infty\,;
\end{equation}
  $\bullet\;$ Pareto and Pareto-like distribution with $K, \alpha>0$:
\begin{equation}
  \oF(x) \sim K\, x^{-\alpha} =:\zeta(x)\,,\;\; x\to\infty\,;
\end{equation}

\item[] \hspace{-7mm} Weibull MDA distributions:\\
  $\bullet\;$ Beta distribution with $0<x<1$, $p,q>0$, $x_E=1$ and $c_{p,q}=\Gamma(p+q)/(\Gamma(p)\Gamma(q))$:
\begin{equation}
\oF(x) = c_{p,q} \int_x^{\infty} \left( u^{p-1}\,(1-u)^{q-1}\right)\mathrm{d}u \sim \frac{c_{p,q}}{q} (1-x)^{q} =: \zeta(x) \,,\;\; x\uparrow x_E\,;
\end{equation}
$\bullet\;$ Pareto-like distribution with finite right endpoint $x_E<\infty$ with $x< x_E$ and $K, \alpha>0$:
   \begin{equation}\label{eq:specdist:E}
   \oF(x)\sim K\, (x_E-x)^{\alpha} =:\zeta(x)\,,\;\; x\uparrow x_E\,.
   \end{equation}
\end{itemize}

For distributions where no explicit representation of the distribution function~$F$ exists, e.g. Gaussian, Lognormal, Gamma, Loggamma, or Beta distributions, the statements $\oF(x)\sim \zeta(x)$ as $x\uparrow x_E$, can be proven by l'H\^opital's rule by using the representation of density~$f=-\oF'$.

Applying Theorem~\ref{Prop:zeta} for these asymptotically equivalent functions~$\zeta$, we now obtain the corresponding auxiliary functions. Note that by selecting an auxiliary function valid for an MDA distribution there is a substantial difference between \vMR and \VR cases: whereas there is no much freedom in \vMR case,
the entire asymptotic equivalence set $\{\psi\sim\psi_{\mathrm{u}}\}$ is available for \VR according to the characterization in Theorem~\ref{Th:psi}(i).
This implies that it is possible to simplify \VRp-valid auxiliary functions to some extent.
We present some convenient auxiliary functions $\psi_{\mathrm{u}}$ and~$\psi$ for selected distributions using the results from Eqs.~\eqref{eq:specdist:A}--\eqref{eq:specdist:E} in the next corollary.

\begin{corollary}\label{Cor:special-distr}
The universal auxiliary functions $\psi_{\mathrm{u}}$ valid for both \VR and \vMRp,
as well as some parsimonious auxiliary functions~$\psi$ valid for \VRs only:
\begin{center}
\begin{tabular}{clclc}
\hline
&  MDA distribution & $\dps\psi_{\mathrm{u}}(x)$ & & $\dps\psi(x)$\\
\hline\\[-2mm]
  \parbox[t]{2mm}{\multirow{7}{*}{\rotatebox[origin=c]{90}{$\,$ Gumbel $\quad$}}} &Gaussian distribution:  & $ x/(1+x^2)$ & & $  1/x$\\[2mm]
  &Lognormal distribution:  & $  x\, \log x /(1+(\log x)^2) $ & & $   x/\log x $\\[2mm]
  &Exponential-like distribution:  & $  1/\lambda $ & &  $  1/\lambda $\\[2mm]
  &Gamma distribution:  & $   x/(\beta x -\alpha +1) $ & & $   1/\beta $\\[2mm]
  &Weibull-like distribution: & $ x/(\tau (\beta x)^{\tau} -\alpha)  $ & & $  x^{1-\tau}/(\tau \beta^{\tau})  $\\[2mm]
  \hline
 & &&&\\[-2mm]
  \parbox[t]{2mm}{\multirow{3}{*}{\rotatebox[origin=c]{90}{Fr\'{e}chet $\;\;\;\,$}}} &Loggamma distribution: & $  x\,\log x/(\alpha \log x\, -\beta +1)  $ & & $  x/\alpha  $\\[2mm]
  &Cauchy distribution: & $  x  $ & & $  x  $\\[2mm]
  &Pareto-like distribution &&&\\
  &($x_E=\infty$): & $   x/\alpha  $ & & $   x/\alpha  $\\[2mm]
  \hline
  & &&&\\[-2mm]
  \parbox[t]{2mm}{\multirow{3}{*}{\rotatebox[origin=c]{90}{Weibull $\;$}}} &Beta distribution: & $  (1-x)/q $ & & $  (1-x)/q $\\[2mm]
 &Pareto-like distribution &&&\\
 &($x_E<\infty$): $\qquad$ & $  (x_E-x)/\alpha  $ & & $  (x_E-x)/\alpha  $\\[2mm]
  \hline
\end{tabular}
\end{center}
\end{corollary}

As we show in Corollary~\ref{Cor:special-distr}, the valid auxiliary functions for Fr\'{e}chet or Weibull MDA distributions depend on the extreme value index~$\gamma$ due to the property $\psi\in\Pg$, which is $\psi(x)\sim \gamma x$ for $x\to\infty$ or $\psi(x)\sim -\gamma (x_E-x)$ for $x\uparrow x_E$, respectively: it holds $\gamma=1/\alpha$ for Loggamma and Pareto-like, $\gamma=1$ for Cauchy, $\gamma=-1/q$ for Beta distributions.

Corollary~\ref{Cor:special-distr} illustrates two major advantages of our idea to construct valid auxiliary functions $\psi_{\mathrm{u}}=-\zeta/\zeta'$ as in Theorem~\ref{Prop:zeta} instead of using $\psi=\oF/F'$ from~\eqref{eq:psi:Gum:3}: first, for many distributions there exists no explicit expression for $\oF(x)/F'(x)$, second, even when it exists, we can often find analytically more simple auxiliary functions based on~$\zeta$.

\begin{remark}\label{Rem:AMF:Criterion} $\quad$\\
(i)$\;\;$ The characterization of the class~$\AMF$ of valid auxiliary function for \vMR in Theorem~\ref{Th:psi}(ii) can be applied as a criterion to check whether a candidate $\psi$ is a \vMRp-valid auxiliary function or not.
This result helps us e.g. to prove that $\psi(x)=\gamma x$ from~\eqref{eq:psi:Fre} is not \vMRp-valid for Loggamma distribution with $\beta\neq1$ and $\gamma=1/\alpha$ in a much more convenient way compared to using the definition as in Remark~\ref{Rem:Examples:Sect3}(ii):
for the universal auxiliary function $\psi_{\mathrm{u}}(x)=x\,\log x/(\alpha \log x\, -\beta +1)$ from Corollary~\ref{Cor:special-distr} we obtain for all $z\in(0,\infty)$ that:
\begin{eqnarray*}
  \left|\int_z^{x}\frac{\psi_{\mathrm{u}}(v) - \psi(v)}{\psi_{\mathrm{u}}(v)\psi(v)}\mathrm{d}v\right| &=& \left|(\beta-1)\,\int_z^{x} \frac{1}{v \log v}\mathrm{d}v\right|\\[2mm]
   &=&\left|\beta-1\right|\, \left(\log(\log x)-\log(\log z)\right)\,\to\, \infty\;\; \text{ for } x\to\infty\,.
\end{eqnarray*}
Consequently, with Theorem~\ref{Th:psi}(ii), $\psi(x)=\gamma x\notin \AMF$, so that $\psi(x)=\gamma x$ is not a \vMRp-valid auxiliary function of the generic Loggamma distribution.\\
(ii)$\;\;$ This convenient criterion based on Theorem~\ref{Th:psi}(ii) can also be applied for various other MDA distributions to prove that the analytical simple choice of~$\psi$ from \VR is \emph{not} a valid choice for \vMRp. For all those MDA distributions from Corollary~\ref{Cor:special-distr} where $\psi$ differs from universal~$\psi_{\mathrm{u}}$ we can easily show that this~$\psi$ is only \VRp-valid but not \vMRp-valid as it holds for $x\to\infty$ that $\int_z^{x} (\psi_{\mathrm{u}}(v) - \psi(v))/(\psi_{\mathrm{u}}(v)\psi(v))\mathrm{d}v =: J(x)-J(z)$ does not have a finite limit~$K_z$ for any fixed $z\in(0,\infty)$. The function~$J$ depends on the particular distribution and is given as $J(x)=-\log x$ for Gaussian, $J(x)=-\log(\log x)$ for Lognormal, $J(x)=(\alpha-1) \log x$ for Gamma, $J(x)=\alpha \log x$ for Weibull-like or $J(x)=(\beta-1) \log(\log x)$ for Loggamma distribution, respectively. As $|J(x)|\to\infty$ as $x\to\infty$ for all these distributions, the corresponding functions~$\psi$ does not belong to $\AMF$ according to Theorem~\ref{Th:psi}(ii).
$\phantom{abc}$ $\hfill \diamond$
\end{remark}

With our results it is possible to identify immediately -- given the functional form of~$\psi$ -- for which family of distributions this~$\psi$ is a valid auxiliary function.
This is a useful supplement to the results of \citet{AFGS08} and \citet{FS12} for estimation
of auxiliary functions~$\psi$ for Gumbel MDA distributions.
In particular, they elaborate on semiparametric estimation methods and assume $\psi$ to be of the following analytical form (see \citealp[Eq.~14]{AFGS08}, \citealp[Eq.~37]{FS12}):
\begin{equation}\label{eq:estimation:psi}
\psi(x)=\frac{1}{c \beta}x^{1-\beta}
\end{equation}
with positive constants $c$ and $\beta$. Hence, estimation of function~$\psi$ is focused on estimation of constants~$c$ and $\beta$; there are various semiparametric estimators proposed in the literature for this purpose.

Our results allow to decide which distributions are covered by this estimation procedure. The table in Corollary~\ref{Cor:special-distr} shows that the functional form~\eqref{eq:estimation:psi} of~$\psi$ as power function is suitable for estimation of the \VRp-valid auxiliary functions for many important Gumbel MDA distributions, a remarkable exception is the Lognormal distribution where $\psi$ is not asymptotically equivalent to~\eqref{eq:estimation:psi}. But this table also shows that the universal auxiliary functions~$\psi_{\mathrm{u}}$ often are not covered by analytical form~\eqref{eq:estimation:psi}.
Hence, our results (see also Theorem 4.1(ii)) show that for most of these Gumbel MDA distributions the functional form as in~\eqref{eq:estimation:psi} is not suitable to be applied for \vMRp-valid auxiliary functions. In order to extend this estimation method also for \vMRp-valid auxiliary functions, the functional form of~$\psi$ should be chosen in a more flexible way.

\section{Conclusions}\label{sec:6}

Both variation representation (\VRp) and von Mises representation (\vMRp) offer a convenient way
to gain limit results for the family of MDA distributions which are essential in extreme value theory.
For this purpose one needs to find so-called auxiliary functions~$\psi$
 valid for \VR and/or for \vMR which are in focus of this paper.
Our major results are as follows: In Theorem~\ref{Th:MisesimpliesVar} we prove that the set of valid auxiliary functions for \vMR is a strict subset of those for \VRp. Moreover, in Theorem~\ref{Prop:zeta} we show how to find valid universal auxiliary functions~$\psi_{\mathrm{u}}$ for the entire MDA distribution family. This result offers us, first, to show in Proposition~\ref{Prop:psi0:2} under which additional assumptions the reciprocal hazard function defines a valid universal auxiliary function, second, to obtain in Proposition~\ref{Prop:psi0} explicit expressions for~$\psi_{\mathrm{u}}$, and, third, to obtain parsimonious auxiliary functions for important distributions in Corollaries~\ref{Cor:EV} and~\ref{Cor:special-distr}. All this allows us to find a universal auxiliary function~$\psi_{\mathrm{u}}$ for each arbitrary MDA distribution.
In Theorem~\ref{Th:psi} we answer the question how much we can deviate from~$\psi_{\mathrm{u}}$ to obtain still a valid auxiliary function.
Here we state ``if and only if''-characterizations for the classes of valid auxiliary functions each for either \VR or \vMRp.

\section{Proofs}\label{sec:7}

\subsection*{Proof of Theorem~\ref{Th:MisesimpliesVar}.}

Let $\psi$ be a valid auxiliary function for \vMR of distribution~$F$, i.e. $\psi$ satisfies~\eqref{eq:repMises} and $\psi\in\Pg$. The proof techniques differ for the three MDA classes and, therefore, the proof is given separately for each of the cases $\gamma=0$, $\gamma>0$, and $\gamma<0$.\\
Firstly, we analyze the case that $F$ is in the Gumbel MDA, i.e. $\gamma=0$. The property $\psi\in\P_0$ guarantees that $\psi$ is absolutely continuous on $(x^\star,\infty)$ with $\lim_{x\uparrow x_E} \psi'(x)=0$. Let $v>0$; by the mean value theorem it follows that there exists a $\xi(x)\in(x,x+v\psi(x))$ such that for $x>x^\star$ it holds that:
\begin{equation}\label{eq:proof:psi:Beurling}
  \frac{\psi(x+v\psi(x))}{\psi(x)}-1=\frac{v\,(\psi(x+v\psi(x))-\psi(x))}{v\,\psi(x)}=v\, \psi'(\xi(x)) \to 0\,,\;
\end{equation}
locally uniformly in $v\in\mathbb{R}$ for $x\to\infty$.
Moreover, $\psi\in\P_0$ implies that
\begin{equation}\label{eq:proof:psioverx:0}
x+z\psi(x) \to x_E \text{ for } x\uparrow x_E
\end{equation}
for arbitrary $z\in\mathbb{R}$, which can be proven as follows. In case of a finite right endpoint~$x_E<\infty$ of distribution $F$, from \vMR in~\eqref{eq:repMises} it follows that $1/\psi(x)\to\infty$ and, hence, $\psi(x)\to0$ for $x\uparrow x_E$. In case of $x_E=\infty$, we have for $x>x^\star$
that
\begin{equation}\label{eq:proof:psioverx}
\frac{\int_{x^\star}^x \psi'(u)\mathrm{d}u}{x} = \frac{\psi(x)-\psi(x^\star)}{x} \;\sim\; \frac{\psi(x)}{x}\,, \; x\to\infty.
\end{equation}
The quotient on the left hand side of equation~\eqref{eq:proof:psioverx} converges to zero, which follows directly if the integral function is bounded or by l'H\^opital's rule if the integral function increases unboundedly, using that $\lim_{x\to \infty}\psi'(x)=0$ from~\eqref{eq:Pg}. Altogether this gives
\begin{equation}\label{eq:proof:psioverx:2}
\lim_{x\to\infty}\frac{\psi(x)}{x} = 0
\end{equation}
and, hence, $ x+z\psi(x) = x(1+z\psi(x)/x) \to \infty$ as $x\to\infty$.\\
Using the properties~\eqref{eq:proof:psi:Beurling} and~\eqref{eq:proof:psioverx:0}, the \vMR in~\eqref{eq:repMises} gives for $z\in\mathbb{R}$ and $x\uparrow x_E$ that:
\begin{eqnarray*}
  \frac{\oF(x+z\psi(x))}{\oF(x)} &=& h_z(x) \exp\left(-\!\!\!\int_x^{x+z\psi(x)}\!\!\!\!\frac{1}{\psi(u)}\mathrm{d}u\right)\\
   &=& h_z(x) \exp\left(-\int_0^{z}\frac{\psi(x)}{\psi(x+v\psi(x))}\mathrm{d}v\right) \to \exp(-z)\,,
\end{eqnarray*}
after substitution $u=x+v\psi(x)$ and with $h_z(x):=c(x+z\psi(x))/c(x)\to 1$;
which proves the result for $\g=0$.\\

Secondly, we analyze the case that $F$ is in Fr\'{e}chet MDA, i.e. $\g>0$. The property $\psi\in\Pg$ implies that $\psi(x)=k(x) x$ for some function $k$ which is eventually positive and continuous with $\lim_{x\to\infty} k(x)=\gamma$. Hence, for each $\epsilon\in(0,1)$ there exists an $x_\epsilon >0$ such that for all $x>x_\epsilon$ it holds that
$|k(x)/\gamma -1|<\epsilon$, and, consequently for $x>x_\epsilon$ and $z>0$:
\begin{equation}\label{eq:proof:psi:inequal}
\!\!\exp\left(\!-\frac{1}{\gamma(1\!-\!\epsilon)} \!\!\int_{x}^{x\!+\!z k(x) x}\hspace{-10mm} 1/u\;\mathrm{d}u\right) \!< \exp\left(\!-\!\!\!\int_{x}^{x\!+\!z k(x) x} \hspace{-10mm} 1/\psi(u)\;\mathrm{d}u\right) \!< \exp\left(\!-\frac{1}{\gamma(1\!+\!\epsilon)}\!\! \int_{x}^{x\!+\!z k(x) x} \hspace{-10mm} 1/u\;\mathrm{d}u\right)\,.
\end{equation}
For $z\in(-1/\g,0)$, this inequality~\eqref{eq:proof:psi:inequal} holds with reversed inequality signs for $x$ large enough such that $x+z k(x) x >x_\epsilon$.
The upper bound of inequality~\eqref{eq:proof:psi:inequal} converges for $x\to\infty$ as follows:
\begin{equation}
  \!\!\exp\left(\!-\frac{1}{\gamma(1\!+\!\epsilon)} \!\!\int_{x}^{x\!+\!z k(x) x} \hspace{-10mm}1/u\;\mathrm{d}u\right) = \exp\left(\!-\frac{1}{\gamma(1\!+\!\epsilon)}\log(1\!+\!zk(x))\right) \to (1\!+\!\gamma z)^{-1/(\gamma(1\!+\!\epsilon))}\,,
\end{equation}

and the lower bound of \eqref{eq:proof:psi:inequal} converges to $(1+\gamma z)^{-1/(\gamma(1-\epsilon))}$. Since this is valid for all $\epsilon\in(0,1)$, we finally get with~\eqref{eq:repMises} that:
\begin{equation}
  \frac{\oF(x+z\psi(x))}{\oF(x)}=h_z(x) \exp\left(-\!\!\!\int_x^{x+z\psi(x)}\!\!\!\!\!\frac{1}{\psi(u)}\mathrm{d}u\right) \to (1+\gamma z)^{-1/\gamma}\,,
\end{equation}
for $x\to\infty$ and $h_z(x):=c(x+z\psi(x))/c(x)\to 1$ with function~$c$ from~\eqref{eq:repMises}. This proves the result for $\gamma>0$.\\
The proof for $\gamma<0$ can be done analogously to those for $\gamma>0$ using the representation $\psi(x)=k(x)\cdot (x_E-x)$ with an eventually positive and continuous function $k$ satisfying $\lim_{x\uparrow x_E} k(x)=-\gamma$, which follows from the property $\psi\in\Pg$.
$\hfill \qed$

$\quad$\\
In following proofs, in particular in those of Proposition~\ref{Prop:psi0}, derivations for the Gumbel case where $F\in $ MDA$(G_0)$ are substantially more extensive than those for the cases where $F\in $ MDA$(G_\g)$ for $\gamma>0$ or $\gamma<0$;
for this reason we extract some proof steps into the following two lemmata:

\begin{lemma}\label{lemma:Itozero}
For $F\in$ MDA$(G_0)$ it holds that
\begin{equation}\label{eq:proof:Itozero:1}
  \lim_{x\uparrow x_E}\! \Is \!=\! \lim_{x\uparrow x_E}\! \IIs \!=\! \lim_{x\uparrow x_E}\! \IIIs \!=\!0\,.
\end{equation}
\end{lemma}
\textbf{Proof.}
In case of a finite right endpoint $x_E<\infty$, equation~\eqref{eq:proof:Itozero:1} is satisfied as both the integrand $\oF$ goes to zero and the integration interval $(x,x_E)$ shrinks to the point~$x_E$. For $x_E=\infty$ we apply that a survival function of a distribution in the Gumbel MDA decreases faster to zero than any power function (a direct consequence of \vMR and $\lim_{x\to\infty}\psi(x)/x=0$ as proven in~\eqref{eq:proof:psioverx:2}), i.e. for any $q\in(0,\infty)$ there exists an $x_{q}$ such that for all $x>x_{q}$ it holds $\oF(x)\leq x^{-q}$. E.g., by choosing $q=4$ we obtain for $x$ sufficiently large that:
\begin{eqnarray} \label{eq:proof:Itozero}
\int_x^{\infty} \int_w^{\infty} \int_v^{\infty} \oF(u)\mathrm{d} u\, \mathrm{d} v\, \mathrm{d} w&\leq& \int_x^{\infty} \int_w^{\infty} \int_v^{\infty}  u^{-4}\,\mathrm{d} u\, \mathrm{d} v\, \mathrm{d} w = \frac{1}{6 x}\to 0\,.
\end{eqnarray}
Hence, the statement in~\eqref{eq:proof:Itozero:1} is proven. $\hfill \qed$

\begin{lemma}\label{lemma:FII:sameaux}
  Let $F\in$ MDA$(G_0)$. Then it follows that also $1-\I\in$ MDA$(G_0)$ and $1-\II\in$ MDA$(G_0)$.
  Let $\psi_1$ be each a \VRp-valid auxiliary function of $F(x)$, $\psi_2$ of $\I$, and $\psi_3$ of $\II$, respectively. Then it holds that $\psi_1(x)\sim\psi_2(x)\sim\psi_3(x)$ as $x\uparrow x_E$. Moreover, there exists a function~$\psi\in\P_0$ (with $\P_0$ as defined in~\eqref{eq:Pg}), such that all the survival functions, namely $\oF(x)$, $\I$ and $\II$, allow for \VRs with this auxiliary function~$\psi$.
\end{lemma}
\textbf{Proof.}
As $F\in$ MDA$(G_0)$ it holds that $\oF$ allows for \vMR with some auxiliary function~$\psi\in\P_0$. Due to Theorem~\ref{Th:MisesimpliesVar}, this auxiliary function~$\psi$ is also valid for the \VR of~$F$.
By using that $\oF$ allows for \VR with~$\psi$, it follows that:
\begin{equation} \label{eq:proof:var}
  \lim_{x\uparrow x_E}\frac{\int_{x+z\psi(x)}^{x_E}\oF(u)\mathrm{d}u}{\I}
  = \lim_{x\uparrow x_E}(1+z\psi'(x)) \lim_{x\uparrow x_E}\frac{\oF(x+z\psi(x))}{\oF(x)} = 1\cdot \exp(-z)
\end{equation}
where we use property~\eqref{eq:proof:psioverx:0} and l'H\^opital's rule in the first step which can be applied as $\I\to0$, see Lemma~\ref{lemma:Itozero}. Hence, we have proven that also $1-\I\in$ MDA$(G_0)$ and that each auxiliary function~$\psi$ for \VR of $\oF$ can be applied for \VR of $\I$, and analogously it can be shown
that\linebreak $1-\II\in$ MDA$(G_0)$ and that each auxiliary function~$\psi$ for \VR of $\I$ can be applied for \VR of $\II$;
moreover, that this $\psi$ can be chosen to satisfy the property~$\P_0$.

Next, consider the family of survival functions on $[0,\infty)$ defined by $\overline{H}_{x}(v):=\oF(x+v)/\oF(x)$ with $x\in(0,x_E)$.
Let $F$ allow for \VR with both~$\psi_1$ and $\tilde{\psi}$. Then it follows that:
\begin{equation}\label{eq:proof:H:g0}
\lim_{x\uparrow x_E}\overline{H}_{x}(z \psi_1(x)) = \lim_{x\uparrow x_E}\overline{H}_{x}(z \tilde{\psi}(x))= \exp(-z)\,.
\end{equation}
By the Khintchine's convergence theorem (see  \citealp[Th.1.2.3]{LLR83}),
also known as ``convergence to types'' theorem (see \citealp[Prop.0.2]{R87}),
the property in~\eqref{eq:proof:H:g0} implies that $\lim_{x\uparrow x_E}(\psi_1(x)/\tilde{\psi}(x))=1$. Hence, this proves that all auxiliary function valid for \VR have to be asymptotically equivalent to each other. And together with the result above that for \VR of $\oF(x)$, $\I$ and $\II$ we can choose the same auxiliary function, it follows that auxiliary functions valid for \VR of $\oF(x)$, $\I$ or $\II$ have to be all asymptotically equivalent as $x\uparrow x_E$. $\hfill \qed$\\

The proof of Theorem~\ref{Prop:zeta} consists of two parts: The first part, that $\psi$ given as in~\eqref{eq:zeta:Th} is a valid universal auxiliary function under the given assumptions on function~$\zeta$, is given in Section~\ref{sec:3} around Eqs.~\eqref{eq:zeta-psi0}--\eqref{eq:zeta-psi0:2}. The second part, the existence of such a function~$\zeta$ with $\oF(x)\sim \zeta(x)$, $x\uparrow x_E$ and $-\zeta/\zeta'\in\Pg$ follows directly from \citet[Rem.~1.2.8]{HF06}.\\
The proofs of Propositions~\ref{Prop:psi:lit} and of~\ref{Prop:psi0:2} are placed after those of Proposition~\ref{Prop:psi0}.

\subsection*{Proof of Proposition~\ref{Prop:psi0}.}

\noindent We have provided the proof idea around Eqs.~\eqref{eq:zeta-psi0}--~\eqref{eq:zeta-psi0:2}, and have shown there how finding a \vMRp-valid auxiliary function of~$\oF$ can be reduced to the problem of finding a \emph{suitable} function~$\zeta\sim \oF$.\\
Firstly, we analyze the case $\gamma=0$.
As $F\in$~MDA$(G_0)$ it holds that $\oF$ allows for \VR with some auxiliary function~$\psi$. From Lemma~\ref{lemma:FII:sameaux} it also follows that $\I$ allows for \VR with the same auxiliary function~$\psi$, which implies that for all $z\neq 0$ it holds:
\begin{eqnarray}\nonumber
  \lim_{x\uparrow x_E}\frac{\int_{x}^{x+z\psi(x)}\oF(u)\mathrm{d}u}{z \I} &=& \lim_{x\uparrow x_E}\left(\frac{1}{z}- \frac{\int_{x+z\psi(x)}^{x_E}\oF(u)\mathrm{d}u}{z \I}\right)\\ \label{eq:proof:g}
  &=& \frac{1-\exp(-z)}{z}\; =:\; g(z)
\end{eqnarray}
where $\lim_{z\to 0} g(z)= \lim_{z\to 0} \exp(-z) = 1$ with l'H\^opital's rule.\\
Now we apply geometrical considerations and interpret the integral $\int_{x}^{x+b}\oF(u)\mathrm{d}u$ with $b>0$ as the signed area of the region bounded by the graph of~$\oF$ and the $x$-axis over the interval $[x,x+b]$. Then, using that $\oF$ is monotonously non-increasing, it follows for $b\in\mathbb{R}$ that:
\begin{equation}\label{eq:proof:Iineq:1}
  \int_{x}^{x+|b|}\oF(u)\mathrm{d}u \;\;\leq\;\; |b|\, \oF(x) \;\;\leq\;\; - \int_{x}^{x-|b|}\oF(u)\mathrm{d}u\,.
\end{equation}
Next, we set $b:=z\psi(x)$, which gives $|b|=|z|\psi(x)$ for $x$ sufficiently large due to the eventual positiveness of auxiliary function~$\psi$, see~\eqref{eq:repVar}. Hence, we obtain for $z\neq 0$ that
\begin{eqnarray}\nonumber
\!\!\!  g(|z|) = \lim_{x\uparrow x_E}\!\frac{\int_{x}^{x+|z|\psi(x)}\oF(u)\mathrm{d}u}{|z| \I} \leq \liminf_{x\uparrow x_E}\!\frac{\psi(x) \oF(x)}{\I} \hspace{39mm}\\  \label{eq:proof:Iineq:2}
\!\!\!  \hspace{11mm} \leq \limsup_{x\uparrow x_E}\!\frac{\psi(x) \oF(x)}{\I} \leq \lim_{x\uparrow x_E}\!\frac{\int_{x}^{x-|z|\psi(x)}\oF(u)\mathrm{d}u}{-|z| \I} = g(-|z|)\,,
\end{eqnarray}
with function~$g$ from~\eqref{eq:proof:g}, and hence, both bounds $g(|z|)$ and $g(-|z|)$ have the limit value~$1$ for $|z|\to 0$.
As inequality~\eqref{eq:proof:Iineq:2} is valid for all $z\neq0$ it follows
with the squeeze theorem that
$\lim_{x\uparrow x_E} (\psi(x) \oF(x)/\I) = 1$,
which implies that the auxiliary function~$\psi$ of $\oF$ satisfies asymptotically for $x\uparrow x_E$:
\begin{equation}\label{eq:proof:psi-sim}
\psi(x)\sim \frac{\I}{\oF(x)}\,.
\end{equation}

Applying the proof steps~\eqref{eq:proof:g}--\eqref{eq:proof:psi-sim} for $\I$ instead of $\oF$, gives that the \VRp-valid auxiliary function of $\I$ is
asymptotically equivalent to\linebreak $\II/\I$, and further, applying the proof steps for\linebreak $\II$ gives that its \VRp-valid auxiliary function is asymptotically equivalent to $\III/\II$.
Together with Lemma~\ref{lemma:FII:sameaux}, which gives that the auxiliary functions for \VRs of $\oF(x)$, $\I$ and\linebreak $\II$ are all asymptotically equivalent,
we obtain for $x\uparrow x_E$ that:
\begin{equation}\label{eq:proof:IIIsim}
  \frac{\I}{\oF(x)} \sim \frac{\II}{\I} \sim \frac{\III}{\II}\,.
\end{equation}

Applying the asymptotic equivalences from~\eqref{eq:proof:IIIsim}, we obtain for $x\uparrow x_E$ that:
\begin{eqnarray*}
  \frac{\left(\II\right)^3}{\oF(x) \left(\III\right)^2} &\sim& \frac{\II \left(\I\right)^2}{\oF(x)\left(\II\right)^2}\\
  &=&\frac{\left(\I\right)^2}{\oF(x) \II}\to 1\,.
\end{eqnarray*}
Next, we define a twice differentiable function $\zeta$ as follows:
\begin{equation}
  \zeta(x):= \frac{\left(\II\right)^3}{\left(\III\right)^2} \sim \oF(x)\,,
\end{equation}
and with $\psi_{\mathrm{u}}(x) := -\zeta(x)/\zeta'(x)$ we obtain according to~\eqref{eq:zeta-psi0:2} an auxiliary function as:
\begin{equation}
  \psi_{\mathrm{u}}(x) = \frac{\II \III }{3 \I \III -2 \left(\II\right)^2}.
\end{equation}

Calculating first and second derivatives of $\zeta$, which are functions in terms of $\oF(x)$, $\I$, $\II$, and $\III$, and then applying property~\eqref{eq:proof:IIIsim} gives for $x\uparrow x_E$:
\begin{eqnarray}
  \zeta'(x) &\sim& - \frac{\I \left(\II\right)^2}{\left(\III\right)^2},\\
    \zeta''(x) &\sim& \frac{\I \left(\II\right)^3}{\left(\III \right)^3}.
\end{eqnarray}

Together we have:
\begin{eqnarray}\nonumber
  \psi'_{\mathrm{u}}(x) &=& \frac{\mathrm{d}}{\mathrm{d}x} \left(-\frac{\zeta(x)}{\zeta'(x)}\right) = \frac{\zeta(x) \zeta''(x)}{(\zeta'(x))^2} -1\\ \label{eq:proof:psi_der}
  &\sim& \frac{\II\,/\,\I}{\III\,/\,\II} -1  \to 0\,,
\end{eqnarray}
for $x\uparrow x_E$, where property~\eqref{eq:proof:IIIsim} is applied in the last step.
Hence, we have proven that $\psi_{\mathrm{u}}\in\Pg$, which concludes the proof for the Gumbel case $\gamma=0$.\\

Secondly, we analyze the case $\gamma>0$.
It holds that $F\in$~MDA$(G_{\gamma})$, $\gamma>0$, if and only if $\oF(x)$
is regularly varying with variation index~$-1/\gamma$ as $x\to\infty$, see \citet[Th.1.2.1(i)]{HF06}, which implies that $\oF(x)/x$
is regularly varying with variation index $\rho:=-(1+\gamma)/\gamma <-1$ as $x\to\infty$. Next, we apply a central result of the regular variation concept, namely Karamata's Theorem (cf. \citealp[Th. 0.6(a)]{R87}),
which leads here to
\begin{equation}\label{eq:proof:Frechet}
\lim_{x\to\infty}\frac{\oF(x)}{\int_x^{\infty}\oF(u)/u\,\mathrm{d}u} = \lim_{x\to\infty}\frac{x\cdot \oF(x)/x}{\int_x^{\infty}\oF(u)/u\,\mathrm{d}u} = -\rho-1 =1/\gamma\,.
\end{equation}
This implies that
\[\oF(x) \sim \frac{1}{\gamma} \int_x^{\infty}\oF(u)/u\,\mathrm{d}u=:\zeta(x)\,,\]
with the derivative $\zeta'(x) = - \oF(x)/(\gamma x)$.
Hence, we obtain:
\begin{equation*}
\psi_{\mathrm{u}}(x) = -\frac{\zeta(x)}{\zeta'(x)}
= \frac{x\,\int_x^{\infty}\oF(u)/u\, \mathrm{d}u}{\oF(x)}\,.
\end{equation*}
As equation~\eqref{eq:proof:Frechet} implies that $\psi_{\mathrm{u}}\in\Pg$ it concludes the proof for $\g>0$.\\

Thirdly, we analyze the case $\gamma<0$.
From \citet[Prop.~1.13]{R87} it follows that $F\in$~MDA$(G_{\gamma})$, $\gamma<0$, if and only if $x_E < \infty$ and
$\oH(z):=\oF(x_E -1/z)$ is regularly varying with variation index~$1/\gamma$ as $z\to\infty$.
This implies
that $\oH(z)/z$ is regularly varying with variation index
$\rho:=(1-\gamma)/\gamma <-1$ as $z\to\infty$. With Karamata's Theorem (cf. \citealp[Th. 0.6(a)]{R87})
it follows that
\begin{equation}\label{eq:proof:Weibull}
\lim_{z\to\infty} \frac{\oH(z)}{\int_z^{\infty} \oH(v)/v\,\mathrm{d}v}  =\lim_{z\to\infty}\frac{z \oH(z)/z}{\int_z^{\infty} \oH(v)/v\,\mathrm{d}v} = -\rho-1 =-\frac{1}{\gamma}\,.
\end{equation}
Using in equation~\eqref{eq:proof:Weibull} the definition of~$\oH$ and then substituting $u:=x_E -1/v$ implies for $z\to\infty$ that
\begin{equation}\label{eq:proof:oF_Weib}
\oF(x_E -1/z) \sim -\frac{1}{\gamma} \int_z^\infty\oF(x_E -1/v)/v\,\mathrm{d}v = -\frac{1}{\gamma} \int_{x_E-1/z}^{x_E}\oF(u)/(x_E-u)\,\mathrm{d}u\,.
\end{equation}
Hence, with $x=x_E-1/z$ we obtain for $x\uparrow x_E$:
\begin{equation}\label{eq:proof:oF_Weib:2}
\oF(x) \sim -\frac{1}{\gamma} \int_{x}^{x_E}\oF(u)/(x_E-u)\,\mathrm{d}u =: \zeta(x)\,,
\end{equation}
where $\zeta$ is absolutely continuous.
Following our idea of proof according to~\eqref{eq:zeta-psi0:2}, with derivative $\zeta'(x) = \oF(x)/(\gamma (x_E-x))$ we obtain:
\begin{eqnarray}\label{eq:proof:oF_Weib:3}
\psi_{\mathrm{u}}(x) &=& -\zeta(x)/\zeta'(x)
= \frac{(x_E-x)\int_{x}^{x_E}\oF(u)/(x_E-u)\,\mathrm{d}u}{\oF(x)}\,.
\end{eqnarray}
With~\eqref{eq:proof:oF_Weib:2} it follows that $\psi_{\mathrm{u}}(x) \sim -\gamma (x_E-x)$ for $x\uparrow x_E$ and, hence, that $\psi_{\mathrm{u}}\in\Pg$.\\
Altogether we have proven that $\psi_{\mathrm{u}}$ defined in Eqs.~\eqref{eq:psi0:G}--\eqref{eq:psi0:W} are  valid auxiliary functions for \vMR of~$\oF$ in Gumbel, Fr\'{e}chet, and Weibull cases, respectively.
By Theorem~\ref{Th:MisesimpliesVar} it follows that these auxiliary functions~$\psi_{\mathrm{u}}$ are also valid for \VR of~$\oF$. Therefore, the result is proven for all~$\g\in\mathbb{R}$. $\hfill\qed$\\

\subsection*{Proof of Proposition~\ref{Prop:psi:lit}.}

\noindent We start to prove the two results for the case $\g=0$, the first for $\psi$ from~\eqref{eq:psi:Gum:1} and the second for $\psi$  from~\eqref{eq:psi:Gum:2}.
Firstly, we obtain for $x>x^\star$ that
\begin{eqnarray} \nonumber
&&\exp\Big(-\int_{x^\star}^x \frac{\oF(v)}{\int_v^{\infty}\oF(u)\mathrm{d}u} \mathrm{d}v\Big)
=  \exp\Big(-\int_{x^\star}^x \Big(-\frac{\mathrm{d}}{\mathrm{d}v}\log \int_v^{\infty}\oF(u) \mathrm{d}u\Big) \mathrm{d}v\Big)\\ \label{eq:Proof:2b}
&=& \exp\Big(-\Big(-\log \int_x^{\infty}\oF(u) \mathrm{d}u + \log \int_{x^\star}^{\infty}\oF(u) \mathrm{d}u\Big)\Big) = \frac{\int_x^{\infty}\oF(u) \mathrm{d}u}{k({x^\star})}\,,
\end{eqnarray}
where $k({x^\star}):=\int_{x^\star}^{\infty}\oF(u) \mathrm{d}u \in(0,\infty)$.
Hence, we get
\begin{equation}\label{eq:Proof:2c}
\oF(x) = \frac{\oF(x)}{\int_x^{\infty}\oF(u)\mathrm{d}u} \int_x^{\infty}\oF(u)\mathrm{d}u= \frac{\oF(x)}{\int_x^{\infty}\oF(u)\mathrm{d}u} k({x^\star}) \exp\Big(-\int_{x^\star}^x \frac{\oF(v)}{\int_v^{\infty}\oF(u)\mathrm{d}u} \mathrm{d}v\Big)\,.
\end{equation}
Equation~\eqref{eq:Proof:2c} states that $\psi(x)=\int_x^{\infty}\oF(u)\mathrm{d}u/\oF(x)$ from~\eqref{eq:psi:Gum:1} is a valid auxiliary function of \vMR in~\eqref{eq:repMises} of $\oF$ if and only if both $\psi\in\P_0$ and $\lim_{x\to\infty} \psi(x)=1/\lambda$ for some $\lambda\in(0,\infty)$, such that $c(x):=k(x^\star)/\psi(x)\to\lambda k(x^\star)\in(0,\infty)$. The right endpoint of~$F$ has  to be $x_E=\infty$, because for $x_E<\infty$ the auxiliary function of \vMR has to satisfy $\lim_{x\uparrow x_E}\psi(x)=0$ because it holds that $\lim_{x\uparrow x_E}\oF(x)=0$.\\
Secondly, for $\psi(x)={\int_x^{\infty}\int_v^{\infty}\oF(u)\mathrm{d}u\mathrm{d}v}/{\int_x^{\infty}\oF(u)\mathrm{d}u}$ from~\eqref{eq:psi:Gum:2} it holds for $x>x^\star$ that:
\begin{eqnarray} \nonumber
\!\!\!\!\!\!\!\!\!\!&&\!\!\!\exp\Big(-\int_{x^\star}^x \frac{\int_w^{\infty}\oF(u)\mathrm{d}u}{\int_w^{\infty}\int_v^{\infty}\oF(u)\mathrm{d}u\,\mathrm{d}v} \mathrm{d}w\Big)\\ \label{eq:Proof:2d}
\!\!\!\!\!\!\!\!\!\!&=& \!\!\! \exp\Big(-\int_{x^\star}^x \Big(-\frac{\mathrm{d}}{\mathrm{d}w}\log \int_w^{\infty}\int_v^{\infty}\oF(u)\mathrm{d}u\,\mathrm{d}v\Big)\mathrm{d}w\Big)
= \frac{\int_x^{\infty}\int_v^{\infty}\oF(u)\mathrm{d}u\,\mathrm{d}v}{\tilde{k}({x^\star})}\,,
\end{eqnarray}
where $\tilde{k}({x^\star}):=\int_{x^\star}^{\infty}\int_v^{\infty}\oF(u)\mathrm{d}u\,\mathrm{d}v$.
Applying~\eqref{eq:proof:IIIsim} from the proof of Proposition~\ref{Prop:psi0} we obtain that
\begin{equation}\label{eq:Proof:2e}
\oF(x) \sim \left(\frac{\int_x^{\infty}\oF(u)\mathrm{d}u}{\int_x^{\infty}\int_v^{\infty}\oF(u)\mathrm{d}u\mathrm{d}v}\right)^2\cdot \int_x^{\infty}\int_v^{\infty}\oF(u)\mathrm{d}u\,\mathrm{d}v =  \frac{\tilde{k}({x^\star})}{\psi^2(x)} \exp\Big(-\int_{x^\star}^x \frac{1}{\psi(v)} \mathrm{d}v\Big)\,,
\end{equation}
which states that \eqref{eq:Proof:2e} is a valid \vMR if and only if both $\psi\in\P_0$ and $\lim_{x\to\infty} \psi(x)=1/\lambda$ for some $\lambda\in(0,\infty)$. \\
Together, we have the necessary and sufficient condition for $\psi$ from~\eqref{eq:psi:Gum:1} and from~\eqref{eq:psi:Gum:2} to be \vMRp-valid.
Hence, it has to hold that $1/\psi(x)=\lambda+\tilde{h}'(x)$ with derivative $\tilde{h}'$ of some twice differentiable function~$\tilde{h}$ where $\tilde{h}'(x)\to 0$ as $x\to\infty$.
Further, $\psi\in\P_0$ implies that $\psi$ is absolutely continuous on $(x^\star, \infty)$ and the derivative~$\psi'$ vanishes asymptotically, which gives that $\tilde{h}''(x)\to 0$ as $x\to\infty$. This leads to
\[\int_{x^\star}^x \frac{1}{\psi(v)}\mathrm{d}v = \lambda x + h(x)\]
where both derivatives $h'(x) \to 0$ and $h''(x) \to 0$ as $x\to\infty$.
Altogether, this proves the results for $\g=0$.\\

Next, we consider the case $\g>0$. A survival function $\oF$ allows for a \vMR with auxiliary function $\psi(x)=\g x$ if and only if
for $x\to\infty$ holds:
\[\oF(x)=c(x)\exp\big(-\int_{x^\star}^x \frac{1}{\g u}\mathrm{d}u\big)= c(x) \exp \big(\log(x^{-1/\g})+\log((x^\star)^{1/\g})\big)\sim K \, x^{-1/\g}\,,\]
with some $c(x)\to c\in(0,\infty)$, and $K:= c\cdot (x^\star)^{1/\g}$. This proves the result for $\g>0$.\\
For the case $\g<0$ it holds that a survival function $\oF$ allows for \vMR with auxiliary function $\psi(x)=-\g (x_E-x)$, $x_E<\infty$ if and only if for $x\uparrow x_E$:
\begin{eqnarray*}
\oF(x)&=&c(x)\exp(-\int_{x^\star}^x \frac{1}{-\g (x_E-u)}\mathrm{d}u\\
&=& c(x) \exp (\log((x_E-x)^{-1/\g})+\log((x_E-x^\star)^{1/\g})\sim K \, (x_E-x)^{-1/\g}\,,
\end{eqnarray*}
with some $c(x)\to c\in(0,\infty)$, and $K:= c\cdot (x_E-x^\star)^{1/\g}$. This proves the result for $\g<0$.
$\hfill \qed$\\

\subsection*{Proof of Proposition~\ref{Prop:psi0:2}.}

\noindent As we now use the result of Theorem~\ref{Prop:zeta} it remains to prove that $\oF(x)/F'(x)$  satisfies the property $\Pg$ from~\eqref{eq:Pg} under the assumptions given in Proposition~\ref{Prop:psi0:2}.\\
We start with the Gumbel case with $\g=0$.
As $F\in$~MDA$(G_0)$ it holds that $\oF$ allows for a \VR with some auxiliary function~$\psi$.
Because $F'$ is assumed to be positive and non-increasing, the proof steps in Eqs. \eqref{eq:proof:Iineq:1}--\eqref{eq:proof:psi-sim} from the proof of Proposition~\ref{Prop:psi0} above can also be applied with $F'$ instead of~$\oF$, which gives that $\psi(x)\sim \oF(x)/F'(x)$, for $x\uparrow x_E$.
This implies that also $F'$ allows for a \VR with the same auxiliary function~$\psi$, because for all $z\in\mathbb{R}$ it holds that:
\begin{equation*}
  \frac{F'(x+z\psi(x))}{F'(x)}\sim \frac{\psi(x)\oF(x+z\psi(x))}{\psi(x+z\psi(x))\oF(x)} \to \exp(-z)\,,\;\; x\uparrow x_E\,,
\end{equation*}
where we use $F'(x)\sim \oF(x)/\psi(x)$ in the first, and both \eqref{eq:proof:psi:Beurling} and the \VR of~$\oF$ in the second step.\\
Further, using that $-F''$ is positive and non-increasing, and that
$F'(x)=\int_x^{x_E}(-F''(u))\mathrm{d}u$
for $x$ close to $x_E$, the proof steps in Eqs. \eqref{eq:proof:Iineq:1}--\eqref{eq:proof:psi-sim} can also be applied for $-F''$ instead of~$\oF$, which gives that $\psi(x)\sim -F'(x)/F''(x)$, for $x\uparrow x_E$.\\
Altogether, we obtain that:
\begin{equation}\label{eq:proof:psi-sim-derivative}
  \psi(x)\sim \frac{\oF(x)}{F'(x)} \sim \frac{-F'(x)}{F''(x)}\,,\;\; x\uparrow x_E\,.
\end{equation}
Under the assumptions of Proposition~\ref{Prop:psi0:2}, it holds that $\oF/F'$ is absolutely continuous in some left neighborhood of~$x_E$ and \eqref{eq:proof:psi-sim-derivative} implies that the derivative of $\oF/F'$ vanishes asymptotically:
\begin{equation}\label{eq:proof:psi-sim-derivative:2}
\frac{\mathrm{d}}{\mathrm{d} x} \left(\frac{\oF(x)}{F'(x)}\right) = -1 -\frac{-F''(x) \oF(x)}{(F'(x))^2} \to 0 \,,\;\; x\uparrow x_E\,.
\end{equation}
Consequently, we obtain that $\oF(x)/F'(x)\in\P_0$.\\
Next, we analyze the Fr\'{e}chet case with $\g>0$. It holds that $F\in$~MDA$(G_{\gamma})$ with $\g>0$ if and only if $\oF$ is regularly varying with variation index $-1/\gamma<0$ as $x\to\infty$, i.e. it holds that
\[\lim_{x\to\infty}\frac{\oF(\lambda x)}{\oF(x)}=\lambda^{-1/\gamma}\,.\]
Because $F$ is assumed to be absolutely continuous with density~$F'$ where $F'$ is non-increasing, we apply the result of Proposition 0.7(b) in \citet[p.21]{R87}, which gives that:
\[\lim_{x\to\infty}\frac{x\cdot \oF'(x)}{\oF(x)}=-1/\gamma\,,\]
and this implies that $\oF(x)/F'(x)\in\Pg$.\\
Finally, the Weibull case with $\g<0$: here we have that $F\in$~MDA$(G_{\gamma})$ if and only if $\oH(x):=\oF(x_E -1/x)$ is regularly varying with variation index $1/\gamma<0$ as $x\to\infty$. Then, analogously to the case $\gamma>0$ discussed above, we obtain that:
\[\lim_{x\to\infty}\frac{x\cdot \oH'(x)}{\oH(x)}=\lim_{x\to\infty}\frac{\oF'(x_E -1/x)}{x\cdot\oF(x_E -1/x)}=1/\gamma\,.\]
Substitution by $s:=x_E -1/x$ gives that:
\[\lim_{s\uparrow x_E}\frac{(x_E -s)\cdot \oF'(s)}{\oF(s)}=1/\gamma\,,\]
and, hence, $\oF(x)/F'(x)\in\Pg$.
Therefore, the first part of the result is proven for all $\g\in\mathbb{R}$.\\
The second part of the result, i.e. that \vMR is satisfied with $c(x)\equiv c$, follows from Eqs. \eqref{eq:zeta-psi0} and~\eqref{eq:zeta-psi0:1} which are satisfied for $\oF=\zeta$.
$\hfill \qed$\\

\subsection*{Proof of Theorem~\ref{Th:psi}.}

\noindent\textbf{(i)}$\;$ Let $F\in$ MDA($G_{\gamma}$) for some $\gamma\in\mathbb{R}$. Based on distribution~$F$ we define a family of survival functions on $[0,\infty)$ by $\overline{H}_{x}(v):=\oF(x+v)/\oF(x)$ with $x\in(0,x_E)$. The fact that $F$ satisfies \VR in~\eqref{eq:repVar} for some auxiliary function $\psi_i$ implies that
\begin{equation}\label{eq:Hx}
\overline{H}_{x}(z \psi_i(x)) \to (1+\gamma z)^{-1/\gamma} \; \text{ for } x\uparrow x_E\,.
\end{equation}
Extending the idea in~\eqref{eq:proof:H:g0} for arbitrary~$\gamma\in\mathbb{R}$, the Khintchine's convergence theorem (see \citealp[Th.1.2.3]{LLR83})
 facilitates the proof of statement~\textbf{(i)} for both necessary and sufficient conditions: If $F$ satisfies \VR in~\eqref{eq:repVar} with two different auxiliary functions $\psi_1$ and $\psi_2$, then according to this theorem, equation~\eqref{eq:Hx} implies that $\psi_1\sim\psi_2$ for $x\uparrow x_E$. Conversely, if $F$ satisfies~\eqref{eq:repVar} with some auxiliary function~$\psi_1$, then this theorem implies for any function~$\psi_2$ with $\psi_2\sim\psi_1$, that \eqref{eq:Hx} holds for both $\psi_1$ and $\psi_2$. This concludes the proof of statement~\textbf{(i)}.\\
\textbf{(ii)}$\;$
Let $\psi\in\AMF$ and $\psi_{\mathrm{u}}$ a universal auxiliary function for~$F$. From
\[\int_z^{x_E}\, \frac{\psi_{\mathrm{u}}(u)-\psi(u)}{\psi_{\mathrm{u}}(u) \psi(u)}\, \mathrm{d}u=:K_z\in(-\infty,\infty)\,,\]
it follows for $x\uparrow x_E$ that
\begin{eqnarray*}
  \exp\left(-\int_z^{x}\frac{1}{\psi(u)} \mathrm{d}u\right)
  &\sim&\exp(-K_z) \exp\left(-\int_z^{x}\frac{1}{\psi_{\mathrm{u}}(u)} \mathrm{d}u\right)\,.
\end{eqnarray*}
We have that $\psi_{\mathrm{u}}$ is a valid auxiliary function for \vMR in~\eqref{eq:repMises} of $\oF$ with some function~$c_0$ and some boundary value~$x_0^\star<x_E$. The value~$z$ from the definition of $\AMF$ can be chosen arbitrarily close to~$x_E$, hence, we choose $z>x_0^\star$ and obtain:
\begin{equation}
  \oF(x)=c_0(x) \exp\left(-\int_{x_0^\star}^{x}\frac{1}{\psi_{\mathrm{u}}(u)} \mathrm{d}u\right) = c(x) \exp\left(-\int_{z}^{x}\frac{1}{\psi(u)} \mathrm{d}u\right)\,,
\end{equation}
where $c(x):=a\cdot c_0(x) \exp(\int_{z}^{x} (\psi_{\mathrm{u}}(u)-\psi(u))/(\psi_{\mathrm{u}}(u)\psi(u)) \mathrm{d}u)\to a\, c _0\exp(K_z)\in(0,\infty)$ with $a:=  \exp(-\int_{x_0^\star}^{z} 1/\psi_{\mathrm{u}}(u) \mathrm{d}u)\in(0,\infty)$, as $\psi_{\mathrm{u}}\in\Pg$ is continuous and positive. Moreover, we have that $\psi\in\AMF$ implies that $\psi\in\Pg$. Consequently, all functions~$\psi\in\AMF$ are valid for \vMRs of~$\oF$.\\
Let, conversely, $\psi$ be a valid auxiliary function for \vMR of~$\oF$.
Per assumption, we have that $\psi_{\mathrm{u}}\in\Pg$ is \vMRp-valid, too, such that we have:
\begin{equation}\label{eq:proof:oF:01}
  \oF(x)=c_0(x) \exp\left(-\int_{x_0^\star}^{x}\frac{1}{\psi_{\mathrm{u}}(u)} \mathrm{d}u\right) = c(x) \exp\left(-\int_{x^\star}^{x}\frac{1}{\psi(u)} \mathrm{d}u\right)\,,
\end{equation}
which implies that
\begin{equation}\label{eq:proof:oF:02}
 \frac{c_0(x)}{c(x)} \exp\left(\int_{x^\star}^{x}\frac{1}{\psi(u)} \mathrm{d}u - \int_{x_0^\star}^{x}\frac{1}{\psi_{\mathrm{u}}(u)} \mathrm{d}u\right) \;=\; 1\,.
\end{equation}
We assume $x_0^\star <x^\star$ without loss of generality. Using the limits $\lim_{x\uparrow x_E}c_0(x)= c_0\in (0,\infty)$ and $\lim_{x\uparrow x_E}c(x)= c\in (0,\infty)$, equation~\eqref{eq:proof:oF:02} implies that:
\begin{eqnarray}\nonumber
 &&\lim_{x\uparrow x_E} \int_{x^\star}^{x}\frac{\psi_{\mathrm{u}}(u)-\psi(u)}{\psi_{\mathrm{u}}(u)\psi(x)} \mathrm{d}u - \int_{x_0^\star}^{x^\star}\frac{1}{\psi_{\mathrm{u}}(u)} \mathrm{d}u\\ \label{eq:proof:log}
   &=& \lim_{x\uparrow x_E} \left(\;\int_{x^\star}^{x}\frac{1}{\psi(u)} \mathrm{d}u - \int_{x_0^\star}^{x}\frac{1}{\psi_{\mathrm{u}}(u)} \mathrm{d}u \right)
= \log c-\log c_0\,.
\end{eqnarray}
From $\psi_{\mathrm{u}}\in\Pg$ it follows that $\psi_{\mathrm{u}}$ is continuous on $(x_0^\star, x_E)$ such that $\int_{x_0^\star}^{x^\star}1/\psi_{\mathrm{u}}(u)\mathrm{d}u=:b\in(0,\infty)$. Hence, with equation~\eqref{eq:proof:log} we obtain:
\begin{eqnarray*}
 \int_{x^\star}^{x_E}\frac{\psi_{\mathrm{u}}(u)-\psi(u)}{\psi_{\mathrm{u}}(u)\psi(x)} \mathrm{d}u = b +  \log c - \log c_0 \in (-\infty,\infty)\,,
\end{eqnarray*}
which gives that $\psi\in\AMF$. Hence, the proof of statement~\textbf{(ii)} is given.\\
The property $\AMF \subsetneq \AVF$ follows from the result in Theorem~\ref{Th:MisesimpliesVar} proven above.\\
$\phantom{abc}$ $\hfill\qed$\\

\subsection*{Proof of Corollary~\ref{Cor:psi0}.}

\noindent Firstly, we prove that the statements in Corollary~\ref{Cor:psi0} are satisfied for any universal auxiliary function~$\psi_{\mathrm{u}}$. For Gumbel case with $x_E=\infty$, we have shown this in~\eqref{eq:proof:psioverx:2}; for Gumbel case with $x_E<\infty$, it follows from:
\begin{equation}\label{eq:proof:psioverx:3}
\lim_{d\downarrow 0}\,\frac{\int_{0}^d \psi'_{\mathrm{u}}(u)\mathrm{d}u}{d} = \lim_{d\downarrow 0}\, \frac{\psi_{\mathrm{u}}(x_E -d)}{d} =\lim_{d\downarrow 0}\, -\psi'_{\mathrm{u}}(x_E-d)=0 \,,
\end{equation}
where we applied l'H\^opital's rule and that $\psi_{\mathrm{u}}(x)\to0$, $x\uparrow x_E$ which follows from the \vMR of survival function~$\oF$. For the Fr\'{e}chet and Weibull cases the results in Corollary~\ref{Cor:psi0} can be obtained as direct consequence of property~$\Pg$ in~\eqref{eq:Pg} for any universal auxiliary function.\\
Secondly, together with Theorem~\ref{Th:psi}(i) by applying the property that all valid auxiliary functions~$\psi$ are asymptotically equivalent to any given universal auxiliary function, the statements of Corollary~\ref{Cor:psi0} are proven for all auxiliary functions~$\psi$.\\
$\phantom{abc}$ $\hfill\qed$\\

\section*{Acknowledgments}
I am grateful to Claudia Kl\"uppelberg for fruitful discussions in which she has drawn my attention to the importance of the research question which is in focus of this paper.

\renewcommand{\baselinestretch}{1.2}


\normalsize

\begin{thebibliography}{99}
{\small \linespread{1}

    \bibitem[Abdous et al.(2008)]{AFGS08}
	\textsc{Abdous, B., Foug\`eres, A.-L., Ghoudi, K., Soulier, P.} (2008).
	\newblock Estimation of bivariate excess probabilities for elliptical models.
	\newblock \textit{Bernoulli}
	\newblock \textbf{14}(4), 1065-1088.\\[-8mm]
		
	\bibitem[Balkema and de Haan(1972)]{BH72}
	\textsc{Balkema, A.A., Haan, de, L.} (1972).
	\newblock On R. von Mises' condition for the domain of attraction of $\exp(-\mathrm{e}^{-x})$.
	\newblock \textit{Ann. Math. Statist.} \textbf{43}(4), 1352-1354.\\[-8mm]
	
	\bibitem[Barbe and Seifert(2016)]{BS16}
	\textsc{Barbe, P., Seifert, M.I.} (2016).
	\newblock A conditional limit theorem for a bivariate representation of a univariate random variable and conditional extreme values.
	\newblock \textit{Extremes} \textbf{19}(3), 351-370.\\[-8mm]

    \bibitem[Beirlant et al.(2017)]{BFR17}
    \textsc{Beirlant, J., Fraga Alves, I., Reynkens, T.} (2017).
    \newblock Fitting tails affected by truncation.
    \newblock \textit{Electron. J. Stat.} \textbf{11}, 2026-2065.\\[-8mm]
	
    \bibitem[D\c{e}bicki et al.(2018)]{DFH18}
    \textsc{D\c{e}bicki, K., Farkas, J., Hashorva, E.} (2018).
    Extremes of randomly scaled Gumbel risks.
    \textit{J. Math. Anal. Appl.} \textbf{458}(1), 30-42.\\[-8mm]

  \bibitem[Embrechts et al.(1997)]{EKM97}
  \textsc{Embrechts, P., Kl\"uppelberg, C., Mikosch, T.} (1997).
  \newblock \textit{Modelling extremal events for insurance and finance.}
  \newblock Springer, Berlin.\\[-8mm]

    \bibitem[Engelke et al.(2019)]{EOW19}
	\textsc{Engelke, S., Opitz, T., Wadsworth, J.} (2019).
	\newblock Extremal dependence of random scale constructions.
	\newblock \textit{Extremes} \textbf{22}(4), 623-666.\\[-8mm]
	
    \bibitem[Fasen et al.(2006)]{FKL06}
    \textsc{Fasen V., Kl\"uppelberg C., Lindner A.} (2006).
    \newblock \textit{Extremal behavior of stochastic volatility models.}
    \newblock In: Shiryaev A.N., Grossinho M.R., Oliveira P.E., Esqu\'ivel M.L. (Eds). Stochastic Finance, 107-155. Springer, New York.\\[-8mm]

	\bibitem[Foug\`eres and Soulier(2012)]{FS12}
	\textsc{Foug\`eres, A.-L., Soulier, P.} (2012).
	\newblock Estimation of conditional laws given an extreme component.
	\newblock \textit{Extremes} \textbf{15}(1), 1-34.\\[-8mm]

    \bibitem[Golosnoy et al.(2012)]{GOS12}
    \textsc{Golosnoy, V., Okhrin, I., Schmid, W.} (2012).
    \newblock Statistical surveillance of volatility forecasting models.
    \newblock \textit{J. Financial Econometrics} \textbf{10}, 513-543.\\[-8mm]
	
	\bibitem[de Haan and Ferreira(2006)]{HF06}
	\textsc{Haan, de, L., Ferreira, A.} (2006).
	\newblock \textit{Extreme value theory. An introduction.}
	\newblock Springer, New York.\\[-8mm]
	
	\bibitem[Hashorva(2012)]{H12}
	\textsc{Hashorva, E.} (2012).
	\newblock Exact tail asymptotics in bivariate scale mixture models.
	\newblock \textit{Extremes} \textbf{15}(1), 109-128.\\[-8mm]

    \bibitem[Hashorva(2013)]{H13}
	\textsc{Hashorva, E.} (2013).
	\newblock Minima and maxima of elliptical arrays and spherical processes.
	\newblock \textit{Bernoulli} \textbf{19}(3), 886-904.\\[-8mm]

    \bibitem[Jan{\ss}en(2010)]{J10}
    \textsc{Jan{\ss}en, A.} (2010).
    \newblock Limit laws for power sums and norms of i.i.d. samples.
    \newblock \textit{Probab. Theory Related Fields} \textbf{146}, 515-533.\\[-8mm]

    \bibitem[Kl\"uppelberg and Lindner(2005)]{KL05}
	\textsc{Kl\"uppelberg, C., Lindner, A.} (2005).
	\newblock Extreme value theory for moving average processes with light-tailed innovations.
	\newblock \textit{Bernoulli} \textbf{11}(3), 381-410.\\[-8mm]

    \bibitem[Kl\"uppelberg and Seifert(2019)]{KS19}
	\textsc{Kl\"uppelberg, C., Seifert, M.I.} (2019).
	\newblock Financial risk measures for a network of individual agents holding portfolios of light-tailed objects.
	\newblock \textit{Finance Stoch.} \textbf{23}, 795–826.\\[-8mm]

    \bibitem[Kl\"uppelberg and Seifert(2020)]{KS20}
	\textsc{Kl\"uppelberg, C., Seifert, M.I.} (2020).
	\newblock Explicit results on conditional distributions of generalized exponential mixtures.
	\newblock \textit{J. Appl. Probab.} \textbf{57}, 760-774.\\[-8mm]

	\bibitem[Leadbetter et al.(1983)]{LLR83}
	\textsc{Leadbetter, M.R., Lindgren, G., Rootz\'en, H.} (1983).
	\newblock \textit{Extremes and related properties of random sequences and processes.}
	\newblock Springer, New York.\\[-8mm]

    \bibitem[McNeil et al.(2015)]{MFE05}
	\textsc{McNeil, A.J., Frey, R., Embrechts, P.} (2015).
	\newblock \textit{Quantitative risk management: Concepts, techniques and tools.}
	Princeton University Press, New Jersey.\\[-8mm]

    \bibitem[Nolde and Wadsworth(2022)]{NW21}
	\textsc{Nolde, N., Wadsworth, J.L.} (2022).
	\newblock Linking representations for multivariate extremes via a limit set.
    \newblock \textit{Adv. in Appl. Probab.} \textbf{54}(3), 688-717\\[-8mm]

    \bibitem[Owada(2017)]{O17}
	\textsc{Owada, T.} (2017).
	\newblock Functional central limit theorem for subgraph counting processes.
	\newblock \textit{Electron. J. Probab.} \textbf{22}(17), 1-38.\\[-8mm]

    \bibitem[Pickands(1975)]{P75}
	\textsc{Pickands, J. III} (1975).
	\newblock Statistical inference using extreme order statistics..
	\newblock \textit{Ann. Statist.} \textbf{3}(1), 119-131.\\[-8mm]

	\bibitem[Resnick(1987)]{R87}
	\textsc{Resnick, S.I.} (1987).
	\newblock \textit{Extreme values, regular variation, and point processes.}
	\newblock Springer, New York.\\[-8mm]
		
	\bibitem[Seifert(2014)]{S14}
	\textsc{Seifert, M.I.} (2014).
	\newblock On conditional extreme values of random vectors with polar representation.
	\newblock \textit{Extremes} \textbf{17}(2), 193–219.\\[-8mm]

    \bibitem[Seifert(2016)]{S16}
	\textsc{Seifert, M.I.} (2016).
	\newblock Weakening the independence assumption on polar components: Limit theorems for generalized elliptical distributions.
	\newblock \textit{J. Appl. Probab.} \textbf{53}(1), 130-145.

}	
\end{thebibliography}
\end{document}